\documentclass{amsart}
\usepackage{graphicx}
\usepackage{amsfonts}
\usepackage{amscd}
\usepackage{amssymb}
\usepackage{xypic}
\newtheorem{thm}{Theorem}[section]
\newtheorem{cor}[thm]{Corollary}
\newtheorem{lem}[thm]{Lemma}
\newtheorem{prop}[thm]{Proposition}
\theoremstyle{definition}
\newtheorem{defn}[thm]{Definition}
\theoremstyle{remark}
\newtheorem{rem}[thm]{Remark}
\numberwithin{equation}{section}
\newcommand{\Matrix}[4]{ \left( \begin{array}{cc}  #1 & #2 \\  #3 & #4 \\ \end{array} \right) }

\newcommand{\Lie}[1]{ {\mathfrak #1} }

\newcommand{\scheme}[1]{\underline{ \mathbf{#1}}}
\newcommand{\variety}[1]{\mathbf{#1}}
\newcommand{\mc}[1]{\mathcal #1}
\newcommand{\Hom}{ {\mbox{Hom}} }
\newcommand{\Ext}{ {\mbox{Ext}} }
\newcommand{\Ind}{ {\mbox{Ind}} }
\newcommand{\nats}{\mathbb N}
\newcommand{\ints}{\mathbb Z}
\newcommand{\rats}{\mathbb Q}
\newcommand{\reals}{\mathbb R}
\newcommand{\complex}{\mathbb C}
\newcommand{\quats}{\mathbb H}
\newcommand{\octs}{\mathbb O}
\newcommand{\hecke}{\mathcal H}
\newcommand{\fourier}{\mathcal F}
\newcommand{\FF}{\mathbb F}
\newcommand{\Gp}{\mathbb G}
\newcommand{\adeles}{\mathbb A}

\newcommand{\Period}{\mathcal{P}}
\newcommand{\isom}{\cong}

\begin{document}

\bibliographystyle{plain}

\title{$D_4$ Modular Forms}%
\author{Martin H. Weissman}%
\date{\today}

\address{Dept. of Mathematics, 970 Evans Hall, University of California, Berkeley, CA 94720-3840}
\email{marty@math.berkeley.edu}%

\thanks{This research was supported by a Liftoff fellowship from the Clay Mathematics Institute and a postdoctoral research fellowship from the National Science Foundation.}%

\begin{abstract}
In this paper, we study modular forms on two simply connected groups of type $D_4$ over $\rats$.  One group, $\variety{G}_s$, is a globally split group of type $D_4$, viewed as the group of isotopies of the split rational octonions.  The other, $\variety{G}_c$, is the isotopy group of the rational (non-split) octonions.  We study automorphic forms on $\variety{G}_s$ in analogy to the work of Gross, Gan, and Savin on $G_2$; namely we study automorphic forms whose component at infinity corresponds to a quaternionic discrete series representation.  We study automorphic forms on $\variety{G}_c$ using Gross's formalism of ``algebraic modular forms''.  Finally, we follow work of Gan, Savin, Gross, Rallis, and others, to study an exceptional theta correspondence connecting modular forms on $\variety{G}_c$ and $\variety{G}_s$.  This can be thought of as an octonionic generalization of the Jacquet-Langlands correspondence.
\end{abstract}

\maketitle

\tableofcontents

\section*{Introduction}

The aim of this paper is to study modular forms on two absolutely simple, simply-connected reductive groups over $\rats$ of type $D_4$.  Originally, I chose to write a paper on this subject because I was interested in the work of Gross, Savin, and Gan (especially in \cite{GGS}) on $G_2$, and I thought I could gain a deep understanding of their methods by writing a paper on the similar phenomena in $D_4$.  This paper was originally meant to be an extended (and hopefully flattering) exercise in imitation.  As a result, the techniques used are for the most part not new, and may be found scattered among the papers on $G_2$.

After writing this paper, I believe that though the techniques may not be new, the results may have independent interest.  Specifically, I believe that the following results in this paper are new, or at least differ from those for $G_2$:
\begin{itemize}
\item
The Fourier coefficients of modular forms on the split group of type $D_4$ are indexed by triples in the ideal class group of imaginary quadratic fields, whose product is the principal class.  Thus, as Siegel's modular forms of degree 2 encode information about ideal class groups of imaginary quadratic fields, $D_4$ modular forms may encode information on the fine structure of these ideal class groups.
\item
The geometry of generalized flag varieties for the split $D_4$ is expressed in terms of octonions.  While essentially contained in the work of J. Tits \cite{Tit}, I hope that the reader finds the simplicity of this approach appealing.
\item
The dual pair correspondence between the split $D_4$ and the $D_4$ which is split at every finite place, and compact at the real place, should be thought of as an octonionic generalization of the Jacquet-Langlands correspondence.  This is a unique phenomenon for $D_4$, and the existence of this correspondence (as suggested to me by B. Gross) is a primary motivation for this paper.
\item
One may construct modular forms on the non-split $D_4$ through the invariant polynomials for the $E_8$ Weyl group.
\end{itemize}

There are many further directions in the study of $D_4$ modular forms.  The action of Hecke operators on the Fourier coefficients of $D_4$ modular forms may be computed as in \cite{GGS} using the relative Satake transform.  More precise information about the $D_4$-$D_4$ dual pair correspondence would be nice, especially relating to values of L-functions.  One may also study modular forms on outer forms of $D_4$ associated to totally real cubic \'etale algebras over $\rats$, and we hope to return to this topic later.

We wish to thank B. Gross, G. Savin, and W.-T. Gan for paving the way for this paper with their research on modular forms for $G_2$, and also their personal advice during the preparation of this paper.  Also, we thank Manjul Bhargava, Stephen DeBacker, Bjorn Poonen, and H.-Y. Loke for answering some questions along the way.  We finally thank the referee for numerous comments and corrections.

\section*{Background and notation}

$\adeles$ will always denote the adeles of the field of rational numbers.  $\hat \ints$ denotes the profinite completion of $\ints$, and is viewed as a $\ints$-algebra, and a subalgebra of $\adeles$.  $\hat \rats = \hat \ints \otimes_\ints \rats$ will denote the finite adeles.

We use different typefaces for schemes, algebraic varieties, and the points of algebraic varieties.  If $R$ is a ring, we use an underlined letter, such as $\scheme{S}$ to denote a scheme over $R$.  If $k$ is a field, we us a boldface letter, such as $\variety{S}$ to denote an algebraic variety over $k$.  Finally, if the field $k$ is fixed, we write $S = \variety{S}(k)$ for the set of $k$-points of $\variety{S}$; if $k$ is a field with a natural topology, such as $\rats_p$ or $\reals$, we endow $S$ with the resulting topology when possible.

If $G$ and $G'$ are groups, and $V,V'$ are representations of $G,G'$ respectively, then we write $V \boxtimes V'$ for the ``external tensor product'' representation of $G \times G'$.  On the other hand, if $V,V'$ are representations of a single group $G$, we write $V \otimes V'$ for the usual tensor product of representations, i.e., $V \otimes V'$ is a representation of $G$, while $V \boxtimes V'$ is a representation of $G \times G$.

If $G$ is a group, and $V$ is a $G$-module, then we write $V^G$ for the subspace of $G$-fixed elements of $V$.  We write $V_G$ for the maximal $G$-invariant quotient of $V$:  $V_G = V / \langle gv - v \rangle$.

We frequently use basic facts about affine algebraic groups, as discussed in Waterhouse \cite{Wat}.  By an algebraic group $\scheme{G}$ over a ring $R$, we always mean an affine algebraic group.  For algebraic groups over a field $k$, we remove the underline, and write $\variety{G}$.  If $R$ is a $k$-algebra, we let $\variety{G}(R)$ denote the $R$-points of $G$.  If $\scheme{G}$ is an affine group scheme over $\ints$, the base change to $\rats$ will be implicitly denoted by removing the underline:  $\variety{G} = \scheme{G} \otimes_{Spec(\ints)} Spec(\rats)$.  $\mu_n$ will denote the group scheme of $n^{\mbox{\small{th}}}$ roots of unity.  $\Gp_a$ will denote the additive group scheme, and $\Gp_m$ will denote the multiplicative group scheme.

If $\variety{G}$ is a reductive algebraic group over $\rats$, we consider three types of automorphic objects for $\variety{G}$.

First, we have the space $\mc{A} = \mc{A}(\variety{G})$ of {\it automorphic forms} on $\variety{G}$, as discussed in the article of Borel \cite{Bor}.  $\mc{A}$ is defined to be the space of smooth functions $f$ on $\variety{G}(\adeles)$ such that:
\begin{itemize}
\item
$f$ is left-invariant under $\variety{G}(\rats)$.
\item
$f$ is right-invariant under some open compact subgroup of $\variety{G}(\hat \rats)$.
\item
$f$ is annihilated by an ideal $J$ of finite codimension in the center of the universal enveloping algebra of the complexified Lie algebra $\Lie{g} \otimes \complex$ of $\variety{G}$.
\item
$f$ is of uniform moderate growth on $\variety{G}(\reals)$.
\end{itemize}
As in the article \cite{GGS}, and in contrast to Borel \cite{Bor}, we do not assume our automorphic forms $f$ to be $K_\infty$-finite for a maximal compact subgroup $K_\infty$ of $G(\reals)$.

Let $\mc{A}^0 = \mc{A}^0(\variety{G})$ denote the cuspidal subspace.  Both $\mc{A}$ and its subspace of cusp forms admit actions of $\variety{G}(\adeles)$ by right translation.

Second, we have the set of {\it automorphic representations} of $\variety{G}$.  We define an {\it automorphic representation} of $\variety{G}$ to be a pair $(\pi, \rho)$, where $\pi$ is an irreducible admissible representation of $\variety{G}(\adeles)$, and $\rho$ (the realization) is a $\variety{G}(\adeles)$-intertwining homomorphism from $\pi$ into the space $\mc{A}(\variety{G})$ of automorphic forms.  We say that $(\pi, \rho)$ is cuspidal if the image of $\rho$ lies in the subspace of cuspidal automorphic forms.  Again, the notion of automorphic representation depends only on the variety $\variety{G}$ and not on the integral structure.

Third, we have the notion of a {\it modular form} for $\variety{G}$.  We define a {\it weight} to be an irreducible smooth representation of $\variety{G}(\reals)$ on a complex Fr\'echet space $W$.  A {\it level} will be an open compact subgroup $\hat K$ of $\scheme{G}(\hat \rats)$; if $0 < N$ is an integer, and $\variety{G}$ comes from a group scheme $\scheme{G}$ over $\ints$ with good reduction everywhere, we associate to $N$ the open compact subgroup
$$\hat K(N) = ker \left( \scheme{G}(\hat \ints) \rightarrow \scheme{G}(\ints / N \ints) \right),$$
and refer to $N$ as a level by abuse of notation.

The space of modular forms of weight $W$ and level $\hat K$ is defined to be the space of $\variety{G}(\reals) \times \hat K$ intertwining homomorphisms from $W \boxtimes \complex$ to the space $\mc{A}(\variety{G})$ of automorphic forms on $\variety{G}$.  The space of cusp forms is defined to be the subspace of homomorphisms whose image lies in the space of cuspidal automorphic forms.

Associated to the open compact subgroup $\hat K$, we may form the global Hecke algebra $\hecke(\hat K)$ consisting of compactly supported, bi-$\hat K$-invariant functions on $\variety{G}(\hat \rats)$.  The space of modular forms of weight $W$ and level $\hat K$ clearly admits an action of this Hecke algebra.  Irreducible Hecke submodules of the space of modular forms yield automorphic representations in the usual manner.

Suppose that $\variety{G}$ is an algebraic group over $\rats$, and $\variety{G}(\rats)$ is discrete in $\variety{G}(\adeles)$ with finite co-volume.  If $f$ is a (measurable) function on $\variety{G}(\adeles)$ that is left-invariant under $\variety{G}(\rats)$, we use the shorthand:
$$\oint_\variety{G} f(g) dg = \int_{\variety{G}(\rats) \backslash \variety{G}(\adeles)} f(g) dg,$$
for the integral with respect to Tamagawa measure.  Note that for unipotent groups $\variety{U}$, Tamagawa measure is normalized so that the compact quotient $\variety{U}(\rats) \backslash \variety{U}(\adeles)$ has volume $1$.

A few such integrals will arise repeatedly, and we mention them here.  If $\variety{G}$ is a reductive group over $\rats$, $\variety{U}$ is a unipotent $\rats$-subgroup of $\variety{G}$, and $f$ is an automorphic form on $\variety{G}$, then we define the $\variety{U}$-constant term of $f$ by the integral:
$$f_\variety{U}(g) = \oint_\variety{U} f(ug) du.$$
More generally, if $\phi$ is a character of $\variety{U}(\adeles)$ which is trivial on $\variety{U}(\rats)$, we define the $\phi$-coefficient of $f$ to be:
$$f_\phi(g) = \oint_\variety{U} f(ug) \overline{\phi(u)} du.$$

If $\variety{U}$ is abelian unipotent, then the characters of $\variety{U}(\adeles)$ trivial on $\variety{U}(\rats)$ can be identified with $\variety{U}(\rats)$, and the Fourier expansion of $f$ reads:
$$f(g) = \sum_{u \in \variety{U}(\rats)} f_{\phi_u}(g).$$

\section{Structure theory}

It is possible to construct a split simply connected group scheme of type $D_4$ over $\ints$ by the methods of Chevalley \cite{Che} using only the root system.  However, we prefer to use a construction that is more special to $D_4$, and which incorporates the most interesting phenomena (triality and octonionic structure) from the start.  From the paper of Gross \cite{Gro}, recall that if $\variety{G}$ is a connected reductive algebraic group over $\rats$, then we say that $\scheme{G}$ is a model for $G$ over $\ints$ if $\scheme{G}$ is a smooth affine group scheme over $\ints$ with general fibre $\variety{G}$, and with good reduction {\it everywhere}.  In \cite{Gro}, Gross gives criteria for the existence of models over $\ints$, and enumerates or classifies these models in many cases.  Specifically, $\variety{G}$ admits a model over $\ints$ if and only if $\variety{G}$ is split over $\rats_p$ for {\it every} prime number $p$.  From \cite{Gro} it follows that:
\begin{itemize}
\item
There is a model $\scheme{G}_s$ of the split simply connected simple group $\variety{G}_s / \rats$ of type $D_4$.  At the real place, it satisfies $\variety{G}_s(\reals) \isom Spin_{4,4}(\reals)$.
\item
There is a model $\scheme{G}_c$ of the simply connected group $\variety{G}_c / \rats$ of type $D_4$ which is split at every finite prime $p$, and which is anisotropic (compact) at the real place, $\variety{G}_c(\reals) \isom Spin_8(\reals)$.
\item
There is a model $\scheme{E}_{8,4}$ of the simply connected group of type $E_8$ which is split at every finite prime $p$, and which is the quaternionic form of real rank $4$ at the real place.
\item
There is a model $\scheme{E}_{7,7}$ of the simply connected group which is split of type $E_7$ at every finite prime $p$, and at the real place as well.
\item
There is a model $\scheme{E}_{7,3}$ of the simply connected group which is split of type $E_7$ at every finite prime $p$, and has real rank $3$ at the real place.
\end{itemize}
These groups are described briefly in \cite{Gro}, and the groups of type $D_n$ are more explicitly constructed in a paper of Goldstine \cite{Gol}.  All of the above integral models are unique in their genus.  Loke explores the dual reductive pair $Spin_{4,4}(\reals) \times_{\mu_2 \times \mu_2} Spin_8(\reals)$ in $E_{8,4}$ over the reals in \cite{Lok}, and in \cite{GJR}, Ginzburg, Jiang, and Rallis study the dual pair $\variety{G}_s \times (\variety{SL}_2^3)$ in $\variety{E}_{7,7}$.  Gross and Savin study dual pairs in $\variety{E}_{7,3}$ in \cite{G-S}.  In this first section we study the structure theory of the groups $\scheme{G}_s, \scheme{G}_c, \scheme{E}_{8,4}, \scheme{E}_{7,7}, \scheme{E}_{7,3}$ in a way which makes these dual pairs more transparent.

\subsection{Octonions}
The connection between exceptional groups, triality, and the octonions can be found in detail in the recent exposition of Baez \cite{Bae}.  The structure of the octonions over $\ints$ was first correctly understood by Coxeter in \cite{Cox}.  An excellent recent exposition on octonions, especially over $\ints$, is provided by the book of Conway and Smith \cite{C-S}.  The structure of the integral split octonions, as well as connections to exceptional groups can be found in the paper of Krutelevich \cite{Kru}.  First, we explicitly describe the groups $\scheme{G}_s, \scheme{G}_c$.  We continue using the subscripts $c$ or $s$ to remind the reader when an object is associated to a group which is compact or split at infinity respectively.

Begin by letting $\quats_c$ denote the division algebra of Hamilton's quaternions, $\quats_c = \reals \oplus \reals i \oplus \reals j \oplus \reals k$, with $ij = -ji = k$, and $i^2 = j^2 = k^2 = -1$.  The main involution on the quaternions is given by:
$$\overline{ (a + bi + cj + dk) } = a - bi - cj - dk.$$
Let $\quats_s$ denote the ``split quaternions'', i.e., the algebra of 2 by 2 matrices with real entries.  The main involution on the matrix algebra is given by:
$$\overline{ \left( \begin{array}{cc}  a & b \\  c & d \\ \end{array} \right) } = \left( \begin{array}{cc}  d & -b \\  -c & a \\ \end{array} \right).$$
Let $Y_s$ denote the maximal order in $\quats_s$ consisting of 2 by 2 matrices with integer entries.  Let $Y_c$ denote the maximal (Eichler) order in $\quats_c$ with $\ints$-basis:
$$1, {1 \over 2} (1+i+j+k), {1 \over 2} (1 + i + j - k), {1 \over 2} (1 - i + j + k).$$
The real numbers $\reals$ can be identified with the center of both $\quats_c$ and $\quats_s$, and both $\quats_c$ and $\quats_s$ split naturally into the direct sum of $\reals$ and $Im(\quats_c)$, $Im(\quats_s)$ respectively.  The space $Im(\quats_s)$ is just the set of trace zero matrices, which can then be identified with the Lie algebra of $SL_2(\reals)$.

We frequently use the ``wildcard'' notation where $\bullet$ may stand for $s$ or $c$.  We may apply the ``Cayley-Dickson process'' to $\quats_\bullet$ to get two eight-dimensional alternative normed algebras with involution over $\reals$.  We let $\octs_\bullet = \quats_\bullet \oplus \quats_\bullet$, with multiplication law:
$$(u,v) \cdot (z,w) = (uz - \bar w v, wu + v \bar z).$$
We define the main involution by:
$$\overline{ (u,v) } = (\bar u, -v).$$
This yields the trace and norm:
\begin{eqnarray*}
Tr(u,v) & = & (u,v) + \overline{(u,v)}  =  (u + \bar u), \\
N(u,v) & = & (u,v) \cdot \overline{(u,v)}.
\end{eqnarray*}
The norm on $\octs_c$ and $\octs_s$ gives a quadratic form on an 8-dimensional real vector space of signature $(8,0)$ and $(4,4)$ respectively.  We call the elements of $\octs_c$ and $\octs_s$ octonions and split octonions respectively.  In $\octs_s$, let $\Omega_s$ be the set of pairs $(u,v)$ with $u$ and $v$ matrices with integer coefficients.  Then $\Omega_s$ is a maximal order in $\octs_s$.  With the symmetric integer-valued bilinear form $\langle \alpha,\beta \rangle = Tr(\bar \alpha \beta)$, $\Omega_s$ is a globally split lattice.  In $\octs_c$, we let $\Omega_c$ denote Coxeter's ring of integral octonions from \cite{Cox}.  Though harder to explicitly describe, $\Omega_c$, endowed with the symmetric integer-valued bilinear form $\langle \alpha,\beta \rangle = Tr(\bar \alpha \beta)$, is isomorphic to the root lattice $E_8$ (after scaling).  We can identify $\quats_\bullet$ as the subalgebra of $\octs_\bullet$ consisting of pairs $(u,0)$.  In this way, $Y_\bullet$ is a subring of $\Omega_\bullet$ as well.

Though multiplication in the algebras $\octs_\bullet$ is not associative, for $\alpha, \beta,\gamma \in \octs_\bullet$ the real number
$$Tr(\alpha \beta \gamma) = Tr(\alpha \cdot (\beta \gamma)) = Tr((\alpha \beta) \cdot \gamma)$$
is well defined.  If moreover, $\alpha,\beta,\gamma \in \Omega_\bullet$, then $Tr(\alpha \beta \gamma) \in \ints$, giving trilinear forms on $\Omega_\bullet$.

\subsection{Integral models}
We consider $(\Omega_\bullet, N)$ as an orthogonal $\ints$-module, from which we get an sequence of group schemes over $\ints$:
$$1 \rightarrow \mu_2 \rightarrow \scheme{Spin}(\Omega_\bullet, N) \rightarrow \scheme{SO}(\Omega_\bullet, N) \rightarrow 1,$$
which is exact in the $fppf$ topology on $Spec(\ints)$.  We refer to the paper of Bass \cite{Bas} for a precise definition of $\scheme{Spin}$ and $\scheme{SO}$ in this case -- we do not simply use the determinant to define $\scheme{SO}$, since that is not the correct notion in the characteristic 2 fibre.
Following Proposition 4.8 of \cite{KPS} we may realize $\scheme{Spin}$ as a subgroup scheme of $\scheme{SO}^3$:
$$\scheme{Spin}(\Omega_\bullet, N) = \{ (\xi, \upsilon, \zeta) \in \scheme{SO}(\Omega_\bullet, N)^3 \colon Tr({}^\xi \alpha {}^\upsilon \beta {}^\zeta \gamma) = Tr(\alpha \beta \gamma) \mbox{ for } \alpha,\beta,\gamma \in \Omega_\bullet \}.$$
We write $\scheme{G}_\bullet$ for $\scheme{Spin}(\Omega_\bullet, N)$, and view points of $\scheme{G}_\bullet$ as triples as above.  Following Gross \cite{Gro}, $\scheme{G}_s$ is the unique model over $\ints$ of the globally split, simply connected simple group of type $D_4$.  $\scheme{G}_c$ is the unique model over $\ints$ of the simply connected simple group of type $D_4$ which is split at every finite place $p$, and which is anisotropic at the real place.

The inclusion of $\mu_2$ in the center of $\scheme{SO}(\Omega_\bullet, N)$ yields an inclusion of $\mu_2^3$ in the center of $\scheme{SO}(\Omega_\bullet, N)^3$.  Define the group scheme:
$$\nu = \{ (\xi,\upsilon,\zeta) \in \mu_2^3 \colon \xi \upsilon \zeta = 1 \}.$$
If $R$ is an integral domain of characteristic zero, then $\nu(R)$ is the finite group $(\ints / 2 \ints)^2$.  Our description of $\scheme{Spin}(\Omega_\bullet, N)$ fixes a natural embedding of $\nu$ in the center of $\scheme{G}_\bullet$.

This construction of $\scheme{G}_\bullet$ can be used to construct the unique model over $\ints$ of the simply connected group of type $E_8$ which is split at every finite place $p$, and which is quaternionic at the real place.  Let $\Lie{g}_\bullet$ denote the Lie algebra over $\ints$ of the group scheme $\scheme{G}_\bullet$, which for all rings $R$ satisfies:
$$\Lie{g}_\bullet \otimes_\ints R \isom ker \left( \scheme{G}_\bullet(R[\epsilon]/\epsilon^2) \rightarrow \scheme{G}_\bullet(R) \right).$$
Let $\Lie{g}_d$ denote the Lie algebra over $\ints$ of the group scheme $(\scheme{G}_c \times_{\nu} \scheme{G}_s)$, where the subscript $\nu$ denotes the natural identification of the copies of $\nu$ in the centers of $\scheme{G}_c$ and $\scheme{G}_s$.  Then $\Lie{g}_d$ contains $\Lie{g}_c \oplus \Lie{g}_s$ with index 4.

The three actions of $\scheme{G}_\bullet$ on $\Omega_\bullet$, which exist by the construction of $\scheme{G}_\bullet$ as a subgroup scheme of $\scheme{SO}(\Omega_\bullet)^3$, yield three infinitesimal actions of the Lie algebra $\Lie{g}_d$ on $\Omega_c \otimes_\ints \Omega_s$.  The algebra structure on $\Omega_c$ and $\Omega_s$ can be further exploited to give a Lie algebra structure on the lattice:
$$\Lie{e}_{8,4} = \Lie{g}_d \oplus \left( \Omega_c \otimes_\ints \Omega_s \right) \oplus \left( \Omega_c \otimes_\ints \Omega_s \right) \oplus \left( \Omega_c \otimes_\ints \Omega_s \right).$$
This construction is completely described by Loke in \cite{Lo2}.

The schematic closure of the subgroup of  $\scheme{GL}(\Lie{e}_{8,4})$ preserving a Killing form and the Lie bracket is an adjoint and simply connected group scheme $\scheme{E}_{8,4}$ over $\ints$, which is split of type $E_8$ over $\rats_p$ at every finite place $p$, and which is the quaternionic real form $E_{8,4}$ over $\reals$.  Essentially by construction, there is an inclusion:
$$\scheme{G}_c \times_{\nu} \scheme{G}_s \hookrightarrow \scheme{E}_{8,4}.$$

Now, we construct the group schemes $\scheme{E}_{7,3}$ and $\scheme{E}_{7,7}$ in a way that makes triality and the dual pair $\scheme{G}_\bullet \times (\scheme{SL}_2^3)$ easy to see.  Note that there is a natural embedding of $\nu$ in the center of $\scheme{SL}_2^3$, so it makes sense to consider the group schemes $\scheme{G}_\bullet \times_\nu \scheme{SL}_2^3$.

Triples $(m_1, m_2, m_3)$ in the Lie algebra $\Lie{sl}_2$ act on triples $(y_1, y_2, y_3) \in Y_s$ by writing:
$$(m_1, m_2, m_3) \cdot (y_1, y_2, y_3) = (m_3 y_1 - y_1 m_2, m_1 y_2 - y_2 m_3, m_2 y_3 - y_3 m_1).$$
This gives a natural action of $Lie(\scheme{G}_\bullet \times_\nu \scheme{SL}_2^3)$ on $(\Omega_\bullet \otimes Y_s)^3$.  Combined with the algebra structures on $\Omega_\bullet \otimes Y_s$, there are natural Lie algebra structures on the lattices:
\begin{eqnarray*}
\Lie{e}_{7,7} & = & Lie(\scheme{G}_s \times_\nu \scheme{SL}_2^3) \oplus \left( \Omega_s \otimes_\ints Y_s \right) \oplus \left( \Omega_s \otimes_\ints Y_s \right) \oplus \left( \Omega_s \otimes_\ints Y_s \right), \\
\Lie{e}_{7,3} & = & Lie(\scheme{G}_c \times_\nu \scheme{SL}_2^3) \oplus \left( \Omega_c \otimes_\ints Y_s \right) \oplus \left( \Omega_c \otimes_\ints Y_s \right) \oplus \left( \Omega_c \otimes_\ints Y_s \right).
\end{eqnarray*}
Again, we point to the work of Loke \cite{Lo2} for a detailed description of these Lie algebras over $\ints$.

The schematic closure of the subgroups of $\scheme{GL}(\Lie{e}_{7,7})$ and $\scheme{GL}(\Lie{e}_{7,3})$ preserving a Killing form and the Lie brackets are simply connected group schemes $\scheme{E}_{7,7}$ and $\scheme{E}_{7,3}$ over $\ints$, which are split of type $E_7$ at every finite place, and have real rank $7$ and $3$ over $\reals$ respectively.  Again, there are dual pair inclusions:
\begin{eqnarray*}
\scheme{G}_s \times_\nu \scheme{SL}_2^3 \hookrightarrow \scheme{E}_{7,7}, \\
\scheme{G}_c \times_\nu \scheme{SL}_2^3 \hookrightarrow \scheme{E}_{7,3}.
\end{eqnarray*}

\section{Modular forms on $\variety{G}_s$}

\subsection{The Heisenberg parabolic}
Choose a maximal torus $\scheme{T}_s$ contained in a Borel subgroup $\scheme{B}_s$ in $\scheme{G}_s$, all over $\ints$.  Let $\Delta_s$ denote the resulting set of simple roots for $\scheme{G}_s$.  The root system of $\scheme{G}_s$ is of type $D_4$ and contains four simple roots
$$\Delta_s = \{ \alpha_0, \alpha_1, \alpha_2, \alpha_3 \},$$
with $\alpha_0$ the ``central'' root, with single edges joining $\alpha_0$ to $\alpha_i$ for $i = 1,2,3$ in the Dynkin diagram.  The highest root $\beta_0$ for $\scheme{G}_s$ can be decomposed:
$$\beta_0 = \alpha_1 + \alpha_2 + \alpha_3 + 2 \alpha_0.$$

Let $\scheme{P}_s$ be the Heisenberg parabolic of $\scheme{G}_s$, associated to the subset of simple roots $\{ \alpha_1, \alpha_2, \alpha_3 \}$.  Let $\scheme{L}_s$ be a Levi component of $\scheme{P}_s$ over $\ints$ containing $\scheme{T}_s$.  Let $\scheme{L}_s'$ denote the derived subgroup of $\scheme{L}_s$.  Then $\scheme{L}'_s$ is isomorphic to $\scheme{SL}_2^3$.  The unipotent radical $\scheme{H}_s$ is a group of Heisenberg type, with center $\scheme{Z}$ one-dimensional.  We call the abelian unipotent 8-dimensional quotient $\scheme{H}_s/\scheme{Z} = \scheme{C}$.  The representation of $\scheme{L}'_s$ by conjugation on $\scheme{C}$ is the tensor cube of the standard $2$-dimensional representation of $\scheme{SL}_2$.  We view points of $\scheme{C}$ as 2 by 2 by 2 cubes, and discuss these in a later section.  For a ring $R$, an element of $\scheme{L}_s(R)$ can be written as a triple of $2 \times 2$ matrices with coefficients in $R$, such as $l = (l_1, l_2, l_3)$, with common non-zero determinant $\det(l)$.

\subsection{Quaternionic discrete series and their continuation}
We work over $\reals$, and so we write $G_s$ for the real Lie group $\variety{G}_s(\reals) \isom Spin_{4,4}(\reals)$.  Let $\varrho$ denote half of the sum of the positive roots, and recall $\beta_0$ is the highest root.

First, we review a few properties of the quaternionic discrete series representations of $G_s$ after Gross and Wallach \cite{G-W} and Wallach \cite{Wal}.  The maximal compact subgroup of $G_s$ is $K \isom SU(2) \times_{\mu_2} (SU(2) \times SU(2) \times SU(2))$.  The representations of $G_s$ that we consider are classified by a lowest $K$-type; these in turn are classified by certain quadruples $(k, \omega_1, \omega_2, \omega_3)$ of non-negative integers.  Let $\complex^2$ denote the standard representation of $SU(2)$ and for $\omega = (\omega_1, \omega_2, \omega_3)$, define the representation of $SU(2) \times SU(2) \times SU(2)$:
$$W_\omega = Sym^{\omega_1}(\complex^2) \boxtimes Sym^{\omega_2}(\complex^2) \boxtimes Sym^{\omega_3}(\complex^2),$$
the external tensor product representation.  We say that a pair $(k, \omega)$ with $\omega = (\omega_1, \omega_2, \omega_3)$ is {\it even} if $k \geq 2$ and $k + \omega_1 + \omega_2 + \omega_3$ is even.  Descending representations from $SU(2) \times SU(2) \times SU(2) \times SU(2)$ to $K$ is described by:
\begin{prop}
There is a bijection between the set of even $(k, \omega)$ and the set of irreducible representations of $K$.  For every even pair $(k,\omega)$, we associate the representation of $K$:
 $$Sym^{k-2}(\complex^2) \boxtimes W_\omega.$$
\end{prop}
When $\omega = 0$, Gross and Wallach describe representations of $G_s$ with lowest $K$-type $Sym^{k-2}(\complex^2) \boxtimes \complex$ in \cite{G-W}; more generally, we summarize some results mentioned in Loke (cf. Theorem 3.3.1 in \cite{Lok}):
\begin{prop}
For $9 \leq k \in \ints$, and $(k, \omega)$ even there is a ``quaternionic'' discrete series representation of $G_s$ with infinitesimal character $\varrho - {k \over 2} \beta_0$, of Gelfand-Kirillov dimension $9$, whose Casselman-Wallach globalization we denote $\pi_{k,\omega}$.  The $K$-finite vectors in $\pi_{k,\omega}$ decompose as a $K$-module via the representations:
$$\bigoplus_{n \geq 0} Sym^{k-2+n}(\complex^2) \boxtimes \left( Sym^n(W_{111}) \otimes W_\omega \right).$$

Even for $2 \leq k < 9$, (with $(k, \omega)$ still even) one may analytically construct representations like the $\pi_{k, \omega}$, which will not however be in the discrete series.  Specifically, if $2 \leq k$, there are smooth representations $\pi_{k, \omega}'$ of $G_s$ on a complex Fr\'echet space, with finite length and infinitesimal character $\varrho - {k \over 2} \beta_0$, whose $K$-finite vectors still decompose according to the previous formula.  This representation $\pi_{k,\omega}'$ may be reducible, but it contains a unique irreducible sub-module $\pi_{k,\omega}$ spanned by the $G_s$-translates of the lowest $K$-type.
\end{prop}
When $\omega=0$, we write $\pi_k$ instead of $\pi_{k,\omega}$.  In particular, the representation $\pi_2$ is unitarizable, and is the minimal representation studied by Kostant in \cite{Kos} of Gelfand-Kirillov dimension $5$.

\subsection{Modular forms}
This section is adapted from Section 7 of Gross, Gan, and Savin \cite{GGS}.  Fix $(k, \omega)$ even, and let $\pi_{k,\omega}$ denote the irreducible representation of $\variety{G}_s(\reals)$ discussed in the last section (a discrete series representation if $9 \leq k$).

Let $\mc{A}_s = \mc{A}(\variety{G}_s)$ denote the space of automorphic forms on $\variety{G}_s$, as discussed in the background section.  $\mc{A}_s^0$ denotes the subspace of cuspidal automorphic forms.
\begin{defn}
Let $\hat K$ denote an open compact subgroup of $\variety{G}_s(\hat \rats)$.  The space of weight $(k,\omega)$ and level $\hat K$ modular forms on $\scheme{G}_s$ is defined to be:
$$\mc{M}_s(k,\omega, \hat K) = \Hom_{\variety{G}_s(\reals) \times \hat K} (\pi_{k,\omega} \boxtimes \complex, \mc{A}_s).$$
When $\omega=0$, we write $\mc{M}_s(k, \hat K)$; we think of $\mc{M}_s(k, \hat K)$ as a space of scalar-valued modular forms of weight $k$, whereas $\mc{M}_s(k,\omega, \hat K)$ is a space of vector-valued modular forms.  The space of weight $(k,\omega)$ and level $1$ modular forms on $\scheme{G}_s$ is defined to be:
$$\mc{M}_s(k,\omega,1) = \Hom_{\variety{G}_s(\reals) \times \scheme{G}_s(\hat \ints)} (\pi_{k,\omega} \boxtimes \complex, \mc{A}_S).$$
The space of weight $(k,\omega)$ cusp forms is defined likewise, replacing $\mc{A}_s$ by $\mc{A}_s^0$, and is denoted $\mc{M}_s^0(k,\omega,\hat K)$.
\label{MoD}
\end{defn}

\subsection{2 by 2 by 2 cubes}
In this section, we follow Bhargava \cite{Bha}, and study the scheme $\scheme{C}$ over $\ints$ which satisfies $\scheme{C}(R) = R^2 \otimes R^2 \otimes R^2$ for any ring $R$.  If $c \in \scheme{C}(R)$, then we think of $c$ as a 2 by 2 by 2 cube $c = (c_{i,j,k})$ of elements of $R$, with $i,j,k \in \{ 0,1 \}$.  If $c$ is a cube, then there are three faces $F_1, F_2, F_3$ of $c$ which contain the entry $c_{0,0,0}$.  Let $F_i'$ denote the faces opposite $F_i$ for $i = 1,2,3$.  We view the $F_i$ and $F_i'$ naturally as 2 by 2 matrices with coefficients in $R$.  From these faces, we define 3 binary quadratic forms:
$$Q_i(x,y) = -\det(F_i x - F_i' y),$$
for $i = 1,2,3$.  The discriminants of the three quadratic forms $Q_i$ are equal to a single $\Delta = \Delta(c)$; thus we call $\Delta(c)$ the discriminant of $c$.  We say that $c$ is non-degenerate if $\Delta(c) \neq 0$.  $\Delta$ is a quartic polynomial map on $\scheme{C}$ with integer coefficients.

The tensor cube of the standard representation yields a natural action of $\scheme{SL}_2^3$ on $\scheme{C}$.  The polynomial $\Delta$ generates the polynomial invariants for this action.  A cube $c \in \scheme{C}(\ints)$ is said to be projective if the three quadratic forms $Q_1, Q_2, Q_3$ are primitive.  In particular, if $c$ is a projective cube, then the greatest common divisor of its entries $c_{i,j,k}$ equals one.   Following the appendix in \cite{Bha}, it is useful to know that projective cubes can be put into a particularly nice form via the action of $\Gamma = \scheme{SL}_2(\ints)^3$:
\begin{prop}
Every projective cube $c \in \scheme{C}(\ints)$ is $\Gamma$-conjugate to some cube $\tilde c$ which satisfies $\tilde c_{0,0,0} = 1$, and $\tilde c_{1,0,0} = \tilde c_{0,1,0} = \tilde c_{0,0,1} = 0$.
\end{prop}
A cube $c$ which already satisfies the conditions in the above proposition is said to be in {\it normal form}.  The discriminant of a cube $c$ in normal form is given by:
$$\Delta(c) = c_{1,1,1}^2 + 4 c_{1,1,0} c_{1,0,1} c_{0,1,1}.$$

For any integer $D \equiv 0,1 (\mbox{mod } 4)$, let $R(D)$ denote the unique quadratic ring over $\ints$ having discriminant $D$ given by:
$$
\begin{cases}
\ints[x] / (x^2) \mbox{ if }  D = 0, \\
\ints + \sqrt{D} (\ints \oplus \ints) \mbox{ if } D \geq 1, \sqrt{D} \in \ints, \\
\ints[(D + \sqrt{D})/2] \mbox{ otherwise }.
\end{cases}
$$

An {\it invertible oriented ideal} in $R(D)$ is a pair $(I, \epsilon)$ where $I$ is a invertible fractional ideal of $R(D)$ in $R(D) \otimes_\ints \rats$, and $\epsilon = \pm 1$.  These form an abelian group by component-wise multiplication.  An invertible {\it principal} oriented ideal is one of the form $((k), sgn(N(k)))$, where $(k)$ is the principal ideal generated by $k \in R(D)$ which is invertible in $R(D) \otimes_\ints \rats$, and $N(k)$ is the norm of $k$ in $\ints$.  The {\it narrow} ideal class group of $R(D)$ is defined to be the quotient group of invertible oriented ideals modulo invertible principal oriented ideals.  When $D > 0$ is square-free, this agrees with the narrow ideal class group of real quadratic fields.  When $D < 0$ is squarefree, the norm form on $R(D)$ is positive-definite, and the narrow ideal class group contains the usual ideal class group with index $2$.

A result essentially known to Gauss is the following:
\begin{prop}
The set of $SL_2(\ints)$ orbits on the set of non-degenerate primitive binary quadratic forms is in bijection with the set of pairs $(D, I)$, where $D \neq 0$ is a discriminant, i.e., a positive integer congruent to $0$ or $1$ mod $4$, and where $I$ is in the narrow ideal class group of $R(D)$.
\end{prop}

A beautiful generalization of the above result is a consequence of Theorem 11 in Bhargava \cite{Bha}:
\begin{prop}
The set of $SL_2(\ints)^3$ orbits on the set of primitive non-degenerate integer cubes of discriminant $D \neq 0$ is in bijection with the set of triples $(I_1, I_2, I_3)$, where $(I_1, I_2, I_3)$ is a triple of narrow ideal classes in $R(D)$ whose product $I_1 \cdot I_2 \cdot I_3$ is the principal class.
\label{Bha11}
\end{prop}
We have seen part of this proposition, since each such integer cube yields three primitive quadratic forms of the same discriminant $D$, and hence three narrow ideal classes in $R(D)$ by the previous result.

\subsection{Heisenberg Whittaker models}

Heisenberg Whittaker models of quaternionic discrete series provide the foundation for the Fourier expansion of modular forms in $\mc{M}_s$.  We review the theory of these Whittaker models here.  As we work only over $\reals$ for the moment, we write $H_s,C,Z, L_s$ for $\variety{H}_s(\reals)$, $\variety{C}(\reals)$, $\variety{Z}(\reals)$, $\variety{L}_s(\reals)$.  Recall that $H_s$ is a Heisenberg group of dimension 9, with abelian quotient $C = H_s / Z$.  The non-degenerate $Z$-valued symplectic form on $H_s$ allows us to identify $C$ with the group of characters $\Hom(H_s, S^1)$ in such a way that $\scheme{C}(\ints)$ is identified with the set of characters which are trivial on $\scheme{H}_s(\ints)$.  Hence, for any $c \in C$, we let $\complex_c$ denote the set of complex numbers, viewed as a representation of $H_s$ via the character $\chi_c$ associated to $c$.  Let $\Delta$ denote the discriminant function on $C$, which is a quartic polynomial.  For any such $c$, and for any even $(k,\omega)$, let $Wh_{k,\omega}(c)$ denote the space of Whittaker models:
$$Wh_{k, \omega}(c) = \Hom_{H_s}(\pi_{k, \omega}, \complex_c).$$

Recall that the Levi component $L_s$ can be identified with the group of triples $l = (l_1, l_2, l_3)$ with $l_i \in GL_2(\reals)$, such that $det(l_1) = det(l_2) = det(l_3)$.  Thus we write $det(l)$ for this common determinant.  The Levi component $L_s$ acts by conjugation on $H_s$, preserving the symplectic form up to scalar, and thus acting on $C$.  As a result, $L_s$ acts on objects related to $C$; if $l \in L_s$ and $c \in C$ then:
\begin{itemize}
\item
We define $l[c] \in C(\reals)$ by $l[c] = l c l^{-1}$.
\item
$\chi_{l[c]}(h) = \chi_c(l^{-1} h l)$ gives an action of $l$ on characters of $H_s$.
\item
If $w \in Wh_{k, \omega}(c)$, and $v$ is in the space of $\pi_{k, \omega}$, we write $l[w](v) = w(l^{-1}(v))$.  Then $l[w] \in Wh_{k, \omega}(l[c])$.
\end{itemize}
In particular, there is a natural action of the stabilizer of $c$ in $L_s$, $L_{s, c}$, on the space $Wh_{k, \omega}(c)$.

Let $W_{k,\omega}$ denote the representation of $L_s$ of highest weight $k \alpha_0 + \omega_1 \alpha_1 + \omega_2 \alpha_2 + \omega_3 \alpha_3$, and $\tilde W_{k, \omega}$ its contragredient.
From Theorem 16 of Wallach \cite{Wal} (and coherent continuation for $k < 9$), we have:
\begin{prop}
Suppose $c \in C$.  If $\Delta(c) < 0$, then $Wh_{k, \omega}(c)$ is isomorphic to $\tilde W_{k,\omega}$ as a representation of the stabilizer $L_{s, c}$ of $c$ in $L_s$.  If $\Delta(c) > 0$, then $Wh_{k, \omega}(c) = 0$.  In particular, if $\omega = 0$, then $Wh_k(c)$ is 1-dimensional if $\Delta(c) < 0$ and the action is $L_{s, c}$ is given by $det^{-3k}$.
\label{Wal16}
\end{prop}
\begin{rem}
It is difficult to ensure that the notion of ``admissible characters'' used by Wallach in \cite{Wal} agrees with the sign of the discriminant used above; such sign errors are very easy to make.  Following advice of W.-T. Gan, we have worked backwards to deduce the correct sign from the existence of a Siegel-Weil type embedding problem.  Namely, we expect the Fourier coefficients of a certain modular form to count embeddings of quadratic rings into the (non-split) octonions with some extra structure.  Such embeddings are only possible if the quadratic rings are imaginary, and thus if the sign of $\Delta$ in the above theorem is negative.
\end{rem}

\subsection{Fourier Expansion}

We consider the Fourier expansion of level 1 modular forms on $\scheme{G}_s$.  Since $\scheme{G}_s$ is a simply-connected simple algebraic group, strong approximation holds, and we may identify:
$$\variety{G}_s(\rats) \backslash \variety{G}_s(\adeles) / \scheme{G}_s(\hat \ints) \isom \scheme{G}_s(\ints) \backslash \variety{G}_s(\reals).$$
For any vector $v \in \pi_{k, \omega}$, and any $f \in \mc{M}_s(k,\omega,1)$, we write $f_v = f(v)$, and view $f_v$ as a function on the single coset space $\scheme{G}_s(\ints) \backslash \variety{G}_s(\reals)$.  Let $\chi_c$ denote any character of the unipotent radical $\variety{H}_s(\reals)$ of the Heisenberg parabolic, which is trivial on $\scheme{H}_s(\ints)$, associated to an element $c$ of $\scheme{C}(\ints)$.  For any continuous function $f$ on $\scheme{H}_s(\ints) \backslash \variety{H}_s(\reals)$, define the (abelian) Fourier coefficient of $f$ by:
$$\fourier_c(f) = \int_{\scheme{H}_s(\ints) \backslash \variety{H}_s(\reals)} f(h) \overline{\chi_c(h)} dh.$$
We may define the linear form $w_c \in Wh_{k,\omega}(c)$ by:
$$w_c(v) = \fourier_c(f_v).$$
Proposition 8.2 of \cite{GGS} carries over to our case to verify that $w_c$ lies in $Wh_{k,\omega}(c)$.  In particular, $w_c = 0$ if $\Delta(c) > 0$.  Moreover if $c' = \gamma[c]$, with $\gamma \in \scheme{L}_s(\ints)$, then we have:
\begin{eqnarray*}
w_{c'}(v) & = & \gamma[w_c](v).
\end{eqnarray*}

When $f$ is a scalar-valued modular form of level 1, $f \in \mc{M}_s(k, 1)$, the uniqueness of Whittaker models in Proposition \ref{Wal16} can be exploited to get scalar Fourier coefficients of $f$ indexed by cubes.  The set of $c \in C(\reals)$ with $\Delta(c) < 0$ forms a single $\variety{L}_s(\reals)$ orbit, and we fix an arbitrary $c_0$ in this orbit.  Since $f$ is scalar-valued, we may also choose an element $w_0 \in Wh_k(c_0)$ which spans $Wh_k(c_0)$.  Hence for any other $c$ of negative discriminant, $w_c$ is a scalar multiple of $l[w_0]$ for some $l \in \variety{L}_s(\reals)$, well-defined up to the stabilizer of $c$ in $\variety{L}_s(\reals)$.  As a result, there are well-defined constants $a_c$ for all $c$ of negative discriminant such that:
$$w_c = a_c \cdot det(l)^{-3k} l[w_0].$$
The factor of $det(l)^{-3k}$ is crucial to make the constant $a_c$ well-defined, even though $l$ is well-defined only up to $L_{s,c}$, using Proposition \ref{Wal16}.

When $f$ is a vector-valued modular form of weight $(k,\omega)$ and of level 1, this may be appropriately generalized:  the coefficients $a_c$ will no longer be scalars, but will have matrix values in $End(W_{k,\omega})$.  When $f$ is scalar-valued and of level 1, we have:

\begin{thm}
Suppose $f \in \mc{M}_s(k, 1)$.  Then the Fourier expansion of $f$ along the Heisenberg parabolic yields constants $a_c$ for every 2 by 2 by 2 cube $c$.  These constants are well-defined, up to a uniform scaling.  The constants $a_c$ vanish if $\Delta(c) > 0$.  For $\Delta(c) < 0$, the constant $a_c$ depends only on the $\Gamma$-orbit of $c$.  Hence, we can associate to each quadruple $(D, I_1, I_2, I_3)$ (a negative discriminant $D = \Delta(c)$ and three oriented ideal classes whose product is the principal class) a constant $a_{(D, I_1, I_2, I_3)}$.
\end{thm}

Let $NCl(D)$ denote the narrow class group of the quadratic ring $R(D)$.  The structure of $NCl(D)$ is completely described by the graph of its group law, i.e., the set of triples $I_1, I_2, I_3$ in $NCl(D)$ whose product is $1$.  Thus modular forms $f \in \mc{M}_s(k,1)$ may provide information on the fine structure of ideal class groups of imaginary quadratic fields.

\section{Modular forms on $\variety{G}_c$}

We now consider modular forms on the group $\variety{G}_c$; recall that $\variety{G}_c(\rats_p)$ is the split simply-connected simple group of type $D_4$ for all $p$, and $\variety{G}_c(\reals)$ is isomorphic to the compact simply-connected simple group $Spin_8(\reals)$.  Thus we may study modular forms on $\scheme{G}_c$ in the framework of Gross's ``algebraic modular forms'' \cite{Gr3}.  With the integral model $\scheme{G}_c$ constructed in the first section, the group of integer points $\scheme{G}_c(\ints)$ has cardinality $2^{14} \cdot 3^5 \cdot 5^2 \cdot 7$, and is isomorphic to the finite group $\Gamma_c := 2^2 \cdot O_8^+(2)$ \cite{Gro}.  The group $\scheme{G}_c(\ints)$ is a central extension of its reduction $\scheme{G}_c(\FF_2)$ by the subgroup $\nu(\ints)$ of order $4$.  More generally, every arithmetic subgroup of $\variety{G}_c(\rats)$ is finite.

\subsection{Algebraic modular forms}
Since $\variety{G}_c$ is an inner form of the split group $\variety{G}_s$ over $\rats$, the irreducible algebraic representations of $\variety{G}_c$ over $\rats$ are parameterized by dominant weights (over an algebraic closure); these are parameterized by non-negative pairs $(k, \omega)$, with $\omega = (\omega_1, \omega_2, \omega_3)$.  Let $V_{k, \omega}$ denote the associated algebraic representation of $\variety{G}_c$ on a rational vector space.

Let $\mc{A}_c = \mc{A}(\variety{G}_c)$ denote the space of automorphic forms on $\variety{G}_c$; for the definition, we refer the reader to the background section.  The condition of moderate growth is unnecessary here, since $\variety{G}_c(\reals)$ is compact.  Fix an open compact subgroup $\hat K$ of $\variety{G}_c(\hat \rats)$; note that $\variety{G}_c(\hat \rats) \isom \variety{G}_s(\hat \rats)$ so the choices of level structure for $\variety{G}_c$ and $\variety{G}_s$ are equivalent.  The space of modular forms of weight $(k,\omega)$ and level $\hat K$ is then:
$$\mc{M}_c(k, \omega, \hat K) = \Hom_{\variety{G}_c(\reals) \times \hat K}(V_{k,\omega} \boxtimes \complex, \mc{A}_c).$$
When $\omega = 0$, we write $\mc{M}_c(k, \hat K)$, and when $\hat K = \scheme{G}_c(\hat \ints)$, we write $\mc{M}_c(k, \omega, 1)$ and view these as modular forms of level $1$.  The complex vector space $\mc{M}_c(k, \omega, \hat K)$ has a natural rational structure, through Gross's theory of ``algebraic modular forms'' in \cite{Gr3}.  We describe this algebraic structure here:

Define the $\rats$ vector space $\mc{M}_c^{alg}(k, \omega, \hat K)$ to be the set of functions $f$ from $\variety{G}_c(\adeles)$ to $V_{k, \omega}$ which are right-invariant by $\hat K$, left-{\it equivariant} by $\variety{G}_c(\rats)$ (with its action on $V_{k,\omega}$), and right-invariant by $\variety{G}_c(\reals)$ too.  Since $\variety{G}_c$ is isomorphic to $\variety{G}_s$ over $\rats_p$ for all $p$, the space of modular forms $\mc{M}_c(k,\omega, \hat K)$ possesses a natural action of the same Hecke algebra $\hecke(\hat K)$ that acts on $\mc{M}_s(k, \omega, \hat K)$.

For any $f \in \mc{M}_c^{alg}(k, \omega, \hat K)$, we define an element $F$ of $\mc{M}_c(k, \omega, \hat K)$ by writing:
$$F_v(g) = \langle g_\infty^{-1} f(g), v \rangle,$$
for all $v \in V_{k, \omega} \otimes \complex$ and $g \in \variety{G}_c(\adeles)$ with archimedean component $g_\infty$.  From Proposition 8.5 of \cite{Gr3}, we have:
\begin{prop}
The map $f \mapsto F$ extends to a $\hecke(\hat K)$-equivariant isomorphism:
$$\mc{M}_c^{alg}(k, \omega, \hat K) \otimes_\rats \complex \isom \mc{M}_c(k, \omega, \hat K).$$
\end{prop}

Thus modular forms on $\scheme{G}_c$ can be studied purely algebraically.  If the $p$ component of $\hat K$ is $\scheme{G}_c(\ints_p)$, then the local spherical Hecke algebra $\hecke(K_p)$ acts on the space of modular forms;  in this setting, it follows from Proposition 8.9 of \cite{Gr3} that every eigenvalue $\lambda$ of a Hecke operator in $\hecke(K_p)$ on $\mc{M}_c(k, \omega, \hat K)$ is algebraic; in fact, it lies in the ring of integers of a CM-field, localized away from $p$.

We study modular forms of level $1$ in more detail.  The double-coset space $\variety{G}_c(\rats) \backslash \variety{G}_c(\adeles) / \variety{G}_c(\reals) \scheme{G}_c(\hat \ints)$ has only one element; this essentially follows from the uniqueness of the $E_8$ lattice as an even unimodular positive-definite lattice of rank $8$.  From Proposition 4.5 of \cite{Gr3}, we have:
\begin{prop}
The space of modular forms of weight $(k,\omega)$ and level $1$ is isomorphic to $V_{k,\omega}^{\Gamma_c}$.
\label{MF1}
\end{prop}
As a trivial first case, there is a one-dimensional $\rats$-vector space of modular forms of level $1$ corresponding to the trivial representation $V = \rats$ of $\variety{G}_c(\rats)$.  But $\Gamma_c$ acts irreducibly on the three fundamental $8$-dimensional representations of $\variety{G}_c(\rats)$, so there are no modular forms of level $1$ of weights $(0,(1,0,0))$, $(0,(0,1,0))$, or $(0,(0,0,1))$.  Of course, there will be modular forms of higher level corresponding to these representations, since the trivial group is an arithmetic subgroup of $\variety{G}_c(\rats)$.

More generally, if $\hat K$ is any level, we fix representatives $g_\delta \in \variety{G}_c(\adeles)$ for the finite collection of double-cosets $\variety{G}_c(\rats) \backslash \variety{G}_c(\adeles) / \variety{G}_c(\reals) \hat K$.  For each $\delta$, we have a finite arithmetic subgroup of $\variety{G}_c(\rats)$:
$$\Gamma_\delta = \variety{G}_c(\rats) \cap g_\delta (\variety{G}_c(\reals) \times \hat K) g_\delta^{-1}.$$
By Proposition 4.5 of \cite{Gr3}, the choice of $g_\delta$ yields an isomorphism, sending $f$ to $\bigoplus f(g_\delta)$:
$$\mc{M}_c^{alg}(k, \omega, \hat K) \isom \bigoplus_\delta V_{k, \omega}^{\Gamma_\delta}.$$

\subsection{Some geometry for $\variety{G}_c$}
We have seen that there are three inequivalent 8-dimensional representations $V_{0,(1,0,0)}, V_{0,(0,1,0)}, V_{0,(0,0,1)}$ of $\variety{G}_c$, all defined over $\rats$.  Let $\octs_\rats = \Omega_c \otimes_\ints \rats$ denote the rational octonion division algebra.  Then we view $\octs_\rats^3$ as the direct sum of the three 8-dimensional representations.

We examine the stabilizers in $\variety{G}_c$ of various special triples in $\octs_\rats^3$.  Consider first the case where all but one of $\alpha, \beta, \gamma$ equals $0$.  Let $\variety{G}_{I}(\alpha)$, $\variety{G}_{II}(\beta)$, $\variety{G}_{III}(\gamma)$ denote the algebraic subgroups of $\variety{G}_c$ stabilizing the triples $(\alpha,0,0), (0,\beta,0), (0,0,\gamma)$ respectively.  All three of these groups become isomorphic to $Spin_7$ over $\reals$.  Thus we call such groups $Spin_7$ subgroups of $\variety{G}_c$ of class $I, II, III$.  By \cite{Var}, there exist precisely three conjugacy classes of subgroups of the real Lie group $Spin_8(\reals)$ isomorphic to $Spin_7(\reals)$ -- the real points of $\variety{G}_{I}(\alpha)$, $\variety{G}_{II}(\beta)$, $\variety{G}_{III}(\gamma)$ represent these three conjugacy classes.

Now, consider the case when only one of $\alpha, \beta, \gamma$ vanishes.  Then the stabilizer of such a vector, e.g. $(\alpha, \beta, 0)$, in $\variety{G}_c$ is a rational algebraic subgroup, which we call $\variety{G}_{I,II}(\alpha, \beta)$, which is an intersection of two $Spin_7$ subgroups of different class.  By Theorem 5 of the third section in \cite{Var}, the intersection of two ``unlike'' $Spin_7(\reals)$ subgroups in $Spin_8(\reals)$ is a subgroup isomorphic to the compact Lie group $G_2$.  Thus we call $\variety{G}_{I,II}(\alpha, \beta)$ a $G_2$ subgroup of $\variety{G}_c$.  We have the following rational version of Theorem 5 of \cite{Var}:
\begin{prop}
Let $\variety{G}_{I}(\alpha)$, $\variety{G}_{II}(\beta)$ denote two $Spin_7$ subgroups of $\variety{G}_c$.  Then there exists a $Spin_7$ subgroup of class $III$, $\variety{G}_{III}(\gamma)$ such that:
$$\variety{G}_{I,II}(\alpha,\beta) = \variety{G}_{I}(\alpha) \cap \variety{G}_{II}(\beta) = \variety{G}_{I}(\alpha) \cap \variety{G}_{II}(\beta) \cap \variety{G}_{III}(\gamma).$$
\end{prop}
\proof
Let $\gamma = \overline{\alpha \beta} = \bar \beta \bar \alpha$ (or any real scalar multiple thereof).  If $g = (\xi, \upsilon, \zeta) \in \variety{G}_{I,II}(\alpha, \beta)$, then $\xi$ stabilizes $\alpha$ and $\upsilon$ stabilizes $\beta$.  To prove the proposition, it suffices to show that $\zeta$ stabilizes $\gamma$; in this case $\variety{G}_{III}(\gamma) \supset \variety{G}_{I,II}(\alpha,\beta)$ and we are done.

To see that $\zeta$ stabilizes $\gamma$, note that $Tr(\alpha \beta \gamma) = Tr(\alpha \beta {}^\zeta \gamma)$ from the definition of $\variety{G}_c$.  We see that $\alpha \beta \gamma = N(\alpha \beta) \in \rats$, and moreover $\alpha \beta ({}^\zeta \gamma - \gamma)$ is totally imaginary since it has zero trace.  Hence the two vectors $\alpha \beta \gamma$ and $\alpha \beta ({}^\zeta \gamma - \gamma)$ are perpendicular in $\octs_\rats$.  Note that $N(\alpha \beta \gamma) = N(\alpha \beta {}^\zeta \gamma)$ and by the Pythagorean Theorem we see that $N(\alpha \beta ({}^\zeta \gamma - \gamma)) = 0$.  Thus ${}^\zeta \gamma - \gamma = 0$.  Hence $\gamma$ is stabilized by $\zeta$ and we are done.
\qed

We have seen that there are essentially three types of $Spin_7$ subgroups of $\variety{G}_c$, which we called class $I, II, III$.  However there is only one type of $G_2$ subgroup of $\variety{G}_c$ by the above proposition -- they arise as intersections of two $Spin_7$ subgroups of different classes, or equivalently as intersections of three $Spin_7$ subgroups which are incident as in the proposition.  Thus we define:
\begin{defn}
Suppose that $(\alpha, \beta, \gamma) \in \octs_\rats^3$ satisfies $(\alpha \beta) \gamma \in \rats$.  Let $\variety{G}_2(\alpha, \beta, \gamma)$ denote the stabilizer of the triple $(\alpha, \beta, \gamma)$ in $\variety{G}_c$.
\end{defn}

Finally, we consider the case of ``generic'' triples $(\alpha, \beta, \gamma)$, in the sense that $(\alpha \beta) \gamma \not \in \rats$.  The stabilizer in $\variety{G}_c$ of such a triple is a subgroup, which by \cite{Var} is isomorphic to $SU(3)$ over $\reals$.  Thus we define:
\begin{defn}
If $\alpha, \beta, \gamma$ are non-zero octonions satisfying $(\alpha \beta) \gamma \not \in \rats$, then let $\variety{SU}_3(\alpha, \beta, \gamma)$ denote the rational algebraic subgroup of $\variety{G}_c$ stabilizing the triple $(\alpha, \beta, \gamma)$.
\end{defn}

\subsection{Periods}
Suppose that $f \in \mc{M}_c^{alg}(k, \omega, \hat K)$ is an algebraic modular form.  Let $\variety{G}'$ be any rational reductive subgroup of $\variety{G}_c$.  The (vector-valued) period of $f$ along $\variety{G}'$ is defined to be:
$$\Period_f^{\variety{G}'}(g) = \oint_{\variety{G}'} (g_\infty')^{-1} f(g'g) dg.$$
We make the following definition:
\begin{defn}
The modular form $f$ is $\variety{G}'$-distinguished if its $\variety{G}'$-period $\Period_f^{\variety{G}'}$ is non-zero.
\end{defn}
As a trivial first case, $f$ is not $\variety{G}_c$-distinguished if and only if $f$ is orthogonal to all constant functions.  Identifying modular forms with vectors in $V_{k, \omega}^\Gamma$ for various finite groups $\Gamma$ yields the following result:
\begin{prop}
Suppose $f$ is a modular form of weight $(k, \omega)$ and level $\hat K$, corresponding to $v \in V_{k, \omega}^{\Gamma_\delta} = f(g_\delta)$.  If $V_{k, \omega}$ has no non-zero vectors invariant under $\variety{G}'(\reals)$, then the period of $f$ must vanish.
\label{Per}
\end{prop}
\proof
The period map for a subgroup $\variety{G}'$ of $\variety{G}_c$ yields a $\variety{G}'(\reals)$-invariant functional on the infinity-type $V_{k, \omega}$ of a modular form.  Hence it is easy to see that if there are no $\variety{G}'(\reals)$-invariant vectors in $V_{k, \omega}$, this functional must vanish.
\qed

Another criterion for $f$ to have vanishing $\variety{G}_2$ periods is given by Corollary 3.1 of \cite{GJR}, which we recall here:
\begin{prop}
Every irreducible admissible generic representation of $\variety{G}_c(\rats_p)$ is not $\variety{G}_2$-distinguished, i.e., has no $\variety{G}_2(\rats_p)$-invariant functionals.
\end{prop}
Hence if $f$ is a modular form on $\variety{G}_c$, which is generic at some finite place $p$, then $f$ is not $\variety{G}_2$-distinguished.

Non-vanishing of periods is more difficult to prove than vanishing.  In the level one case, some non-vanishing results can be proven without too much difficulty.  Consider the triples $(1,0,0)$, $(0,1,0)$, and $(0,0,1)$ of octonions first.  The stabilizers of these triples are $Spin_7$ subgroup schemes $\scheme{G}_I(1,0,0)$, $\scheme{G}_{II}(0,1,0)$, $\scheme{G}_{III}(0,0,1)$ in $\scheme{G}_c$.  The stabilizer of the triple $(1,1,1)$ in $\scheme{G}_c$ is a group scheme $\scheme{G}_2(1,1,1)$.  All of these subgroup schemes have good reduction everywhere -- they are good integral models of the simply connected simple groups of type $B_3$ and $G_2$.

Also, fix a square root of $-1$, $j$, in the integral octonions $\Omega_c$.  The stabilizer of the triple $(1,1,j)$ is an $SU_3$ subgroup scheme $\scheme{SU}_3(1,1,j)$ in $\scheme{G}_c$.  We begin with:
\begin{lem}
The unitary group $\scheme{SU}_3(1,1,j)$ has class number 1, i.e., it satisfies:
$$\# \{ \variety{SU}_3(1,1,j)(\rats) \backslash \variety{SU}_3(1,1,j)(\hat \rats) / \scheme{SU}_3(1,1,j)(\hat \ints) \} = 1.$$
\end{lem}
\proof
We describe the unitary group scheme $\scheme{SU}_3(1,1,j)$ more explicitly.  Let $\FF$ denote the oriented ``Fano plane'':  the finite projective plane with 7 points, and 7 lines.  Fix a basis $\{ 1, e_0, e_1, \ldots, e_6 \}$ of $\octs_\rats$, in which $e_i^2 = -1$ for all $0 \leq i \leq 6$, and where we identify $e_0, \ldots, e_6$ with points $p_0, \ldots, p_6$ of $\FF$.  If $i,j,k \in \octs_\rats$ are square roots of $-1$, we call $(i,j,k)$ a ``quaternion triple'' if they satisfy the familiar relations $ij  = -ji = k$.  The basis of $\octs_\rats$ may be chosen so that $(e_i, e_j, e_k)$ is a quaternion triple in $\octs_\rats$ if $p_i, p_j, p_k$ are oriented collinear points on $\FF$.  This fully describes the multiplication table of $\octs_\rats$ with respect to the chosen basis.

The basis elements $\{ 1, e_0, \ldots, e_6 \}$ are contained in Coxeter's order $\Omega_c$, though they do not span it.  Identifying $j$ with $e_0$ (all square roots of $-1$ are conjugate under $\scheme{G}_c$), let $D$ denote the sublattice of $\Omega_c$ orthogonal to $\ints[e_0]$.  As a lattice, $D$ is isomorphic to the $D_6$ root lattice, and we describe $D$ as:
$$D = \ints\mbox{-span} \left\{ {1 \over 2} (e_i \pm e_j) \right\}, \mbox{ for } 1 \leq i,j \leq 6.$$
$D$ is closed under (left) multiplication by $\ints[e_0]$, and as a $\ints[e_0]$-module it is free of rank three with basis:
$$b_1 = {{e_1 + e_3} \over 2}, b_2 = {{e_2 + e_4} \over 2}, b_3 = {{e_5 - e_6} \over 2}.$$
On $D$ there is a natural $\ints[e_0]$-valued Hermitian form given by:
$$h(d_1, d_2) = - d_1 d_2 + e_0 (d_1 d_2) e_0.$$
An easy calculation shows that $h(b_i, b_j) = \delta_{ij}$, so that $D$ is the simplest possible Hermitian lattice of rank $3$ over $\ints[e_0]$.  The group scheme $\scheme{SU}_3(1,1,j)$ is precisely the unitary group scheme of $D$.  The lemma now follows directly from a result of K. Iyanaga \cite{Iya} (who shows precisely that this unitary group scheme has class number 1).
\qed

In the level 1 case, we can now prove that some periods do not vanish:
\begin{thm}
Suppose that $f$ is an algebraic modular form of level $1$ and weight $(k,\omega)$, corresponding to $v \in V_{k, \omega}^{\Gamma_c}$.  Then if $v$ is invariant under $\variety{G}'(\reals)$ where $\variety{G}'$ is one of $\variety{G}_I(1,0,0)$, $\variety{G}_{II}(0,1,0)$, $\variety{G}_{III}(0,0,1)$, $\variety{G}_2(1,1,1)$, or $\variety{SU}_3(1,1,j)$, then the period $\Period_f^{\variety{G}'}$ is not zero.
\label{PerTh}
\end{thm}
\proof
By the analysis in Proposition 5.5 of \cite{G-S}, to prove our theorem it suffices to show that:
$$\variety{G}'(\adeles) = \variety{G}'(\rats) \variety{G}'(\reals) \scheme{G}'(\hat \ints).$$
For $\variety{G}'$ of type $Spin_7$ or $G_2$, this follows from the uniqueness of integral models discussed in \cite{Gro}.  For $\variety{G}'$ of type $SU_3$, we apply the previous lemma, which shows that the particular unitary group $\variety{SU}_3(1,1,j)$ has class number $1$.
\qed

\subsection{Some branching rules}
In order to determine whether periods of modular forms on $\variety{G}_c$ vanish or not, it suffices by Proposition \ref{Per} and Theorem \ref{PerTh} to understand when representations $V_{k, \omega}$ have fixed vectors when restricted to various reductive subgroups $\variety{G}'$.  Moreover, it suffices to consider these branching problems over the reals.  We consider this problem when $\variety{G}'$ is a $Spin_7$ subgroup, a $G_2$ subgroup, and an $SU_3$ subgroup obtained as before as the stabilier of a triple $(\alpha, \beta, \gamma)$ of octonions.

All groups in this section will be compact simply-connected real Lie groups; we simply write $Spin_8$, $Spin_7$, $G_2$, and $SU_3$ for these groups.  The irreducible representations of $Spin_8$ are the $V_{k,\omega}$ we have already discussed, with $k \geq 0$, and $\omega = (\omega_1, \omega_2, \omega_3)$.  The irreducible representations of $Spin_7$ are indexed by triples $(m, n, r)$ of non-negative integers in the following way:  the triple $(1,0,0)$ corresponds to the 7-dimensional representation, the triple $(0,1,0)$ corresponds to the 21-dimensional adjoint representation, and the triple $(0,0,1)$ corresponds to the 8-dimensional spin representation.  General triples $(m,n,r)$ correspond to representations whose highest weight is the appropriate linear combination of the fundamental weights for the aforementioned representations.  Write $U_{m,n,r}$ for this irreducible representation.

The irreducible representations of $G_2$ correspond to pairs $(p,q)$, where the first coordinate corresponds to the 7-dimensional representation, and the second coordinate corresponds to the 14-dimensional adjoint representation.  Write $T_{p,q}$ for the irreducible representation parametrized as such.

\begin{prop}
Suppose that $\omega = 0 = (0,0,0)$ and $k > 0$.  Then $Res \downarrow_{G_2}(V_{k,\omega})$ does not contain the trivial representation.  However $Res \downarrow_{SU_3}(V_{k,\omega})$ does contain the trivial representation.
\end{prop}
\proof
Choose any $Spin_7$ subgroup of $Spin_8$ containing $G_2$ (there are three conjugacy classes of such $Spin_7$ subgroups, by Varadarajan \cite{Var}).  Then by a well-known branching law, the representation $V_{k, 0}$ restricts as follows:
$$Res \downarrow_{Spin_7} V_{k, 0} = \bigoplus_{m + n = k} U_{m,n,0}.$$
By the branching formula in Theorem 3.4 of \cite{McG}, the representation $Res \downarrow_{G_2} U_{m,n,0}$ contains no trivial representation if $m + n > 0$.  Hence we see that $Res \downarrow_{G_2} V_{k, 0}$ contains no trivial representation.

Now by a branching formula in \cite{Sav}, we see that $Res \downarrow_{SU_3} T_{p,q}$ contains the trivial representation if and only if $q = 0$.  Using McGovern's formula in Theorem 3.4 of \cite{McG} again, it follows that $Res \downarrow_{SU_3} U_{m,n,0}$ always contains the trivial representation.  Hence $Res \downarrow_{SU_3} V_{k, \omega}$ always contains the trivial representation for $\omega = (0,0,0)$.
\qed

If $V_{k, \omega}$ is centrifugal, i.e., $k = 0$, then it follows that $Res \downarrow_{Spin_7} V_{k, \omega}$ contains the trivial representation of $Spin_7$ for some $Spin_7$ subgroup of $Spin_8$.  Hence we have:
\begin{prop}
If $k = 0$, then $Res \downarrow_{G_2} V_{k, \omega}$ contains the trivial representation.
\end{prop}

\begin{cor}
Suppose that $f$ is an algebraic modular form of level $1$ and weight $(k, \omega)$.  If $k = 0$, then $f$ has non-vanishing period along a $G_2$ subgroup.  If $\omega = (0,0,0)$ and $k > 0$, then the periods of $f$ along $G_2$ subgroups vanish, but $f$ has a non-vanishing period along an $SU_3$ subgroup.
\label{PeC}
\end{cor}
\proof
This follows from the last two propositions, together with Theorem \ref{PerTh}.
\qed

\section{Local theta correspondences}
In this section, we study a theta correspondence which can be used to construct modular forms on $\scheme{G}_s$ and $\scheme{SL}_2^3$ from those on $\scheme{G}_c$.  We use the dual pairs:
$$\scheme{G}_s \times_\nu \scheme{G}_c \hookrightarrow \scheme{E}_{8,4},$$
$$\scheme{SL}_2^3 \times_\nu \scheme{G}_c \hookrightarrow \scheme{E}_{7,3}.$$

As $\variety{G}_s$ and $\variety{G}_c$ are inner forms of each other, their Langlands dual groups coincide, and we fix an identification:
$${}^L \variety{G}_s \isom {}^L \variety{G}_c.$$

\subsection{Real Correspondence}

We describe local results over $\reals$ in this section.  In the work of Gross and Wallach \cite{G-W}, they construct the minimal representation $\Pi_\reals$ of the quaternionic group $\variety{E}_{8,4}(\reals)$ by continuation of quaternionic discrete series.  Recall that irreducible representations of $\variety{G}_c(\reals)$ are indexed by pairs $(k, \omega)$, $\omega = (\omega_1, \omega_2, \omega_3)$, with $k, \omega_i$ non-negative.  Let $V_{k,\omega}$ be the irreducible representation associated to the pair $(k, \omega)$.  If $k = 0$, we call $V_{k, \omega}$ a {\it centrifugal} representation -- these will play a special role in the theta correspondence.  Also, recall that for $(k, \omega)$ even, $\pi_{k,\omega}$ is a quaternionic discrete series representation for $k \geq 9$.

The following is one of the main results of the paper of Loke \cite{Lok}:
\begin{thm}
Upon restriction to the dual pair $\variety{G}_s(\reals) \times_{\nu} \variety{G}_c(\reals)$, the representation $\Pi_\reals$ decomposes as a direct sum over the set of irreducible representations of $\variety{G}_c(\reals)$:
$$\mbox{Res }(\Pi_\reals) = \bigoplus_{(k, \omega)} \Theta(V_{k,\omega}) \boxtimes V_{k,\omega},$$
where each $\Theta(V_{k, \omega})$ is isotypic, consisting of the single quaternionic discrete series $\pi_{\vert \omega \vert +2k+10, \omega }$ with finite multiplicity equal to $k + 1$, where $\vert \omega \vert = \omega_1 + \omega_2 + \omega_3$.
\label{LoT}
\end{thm}
In particular, the pairing above is perfect for centrifugal representations, and the trivial representation of $\variety{G}_c(\reals)$ is paired with the quaternionic discrete series of weight $k = 10$, and $\omega = 0$.  Finally, as remarked in 1.7 of \cite{Lok}, this pairing of representations is the same as that predicted by Langlands functoriality.

Now let $\Pi_\reals'$ be the minimal representation of $\variety{E}_{7,3}(\reals)$.  Let $D_k$ denote the holomorphic discrete series representation of $\variety{SL}_2(\reals)$ whose minimal $K$-type corresponds to the positive integer $k$.  A local result due to Gross and Savin, Proposition 3.3 of \cite{G-S}, is:
\begin{thm}
Upon restriction to the dual pair $\variety{SL}_2(\reals)^3 \times_\nu \variety{G}_c(\reals)$, the representation $\Pi_\reals'$ decomposes as a direct sum over the set of irreducible representations $V_{k, \omega}$ of $\variety{G}_c(\reals)$ with $k = 0$:
$$\mbox{Res }(\Pi_\reals') = \bigoplus_{\omega} \Theta'(V_{0, \omega}) \boxtimes V_{0, \omega},$$
where each $\Theta'(V_{0,\omega})$ is given by:
$$\Theta'(V_{0, \omega}) = D_{4 + \vert \omega \vert - \omega_1} \boxtimes D_{4 + \vert \omega \vert - \omega_2} \boxtimes D_{4 + \vert \omega \vert - \omega_3}.$$
\end{thm}
While $\Theta'$ is not functorial in the usual sense, it might be referred to as ``backwards functoriality''.  Note that only centrifugal representations occur in the theta correspondence for $\Pi_\reals'$.

The remainder of this section will be devoted to the $p$-adic versions of the above theorems, in the spherical tempered case.

\subsection{Split $D_4$ geometry}
We begin the process of proving the p-adic theta correspondence,  using methods found in the work of Magaard and Savin \cite{M-S}.  The first step will be geometric.  The geometry of ``amber spaces'' in the 27-dimensional module for $E_6$, discussed in \cite{Asc} and used in \cite{M-S}, will be replaced by a suitable geometry for $D_4$ originated by Tits in \cite{Tit}.

Recall that the algebraic groups $\variety{G}_c$ and $\variety{G}_s$ are isomorphic and split over $\rats_p$.  Thus we write $G$ for the $\rats_p$-points of these groups, working consistently over $\rats_p$ in this section.  Following a tradition of abuse, we call subgroups of $G$ parabolic subgroups, tori, etc..., if they are the $\rats_p$ points of such algebraic subgroups of $\variety{G}$.

Fixing a pinning $T \subset B \subset G$ of $G$, there are 16 standard parabolic subgroups of $G$, corresponding to subsets of the set of simple roots.  The geometric interpretation of these parabolic subgroups, as stabilizers of certain flags, can be expressed in the language of \cite{Tit}.  Of course, viewing $G$ as a classical group, the parabolic subgroups have an interpretation as stabilizers of isotropic flags in (any of) the 8-dimensional algebraic representations of $G$.  It is more canonical to view parabolic subgroups as stabilizers of flags in the full 24-dimensional isotopy representation of $G$, as this does not single out a single 8-dimensional representation.

Thus we prefer, and in fact are required in what comes later, to view $G$ as an exceptional group.  In this section, we write $\octs_p$ for the split (and only) octonion algebra over $\rats_p$.  With this in mind, we define:
\begin{defn}
Let $i,j,k$ be non-negative integers.  Let ${\mathcal Fl}_{i,j,k}$ denote the set of triples $(A,B,C)$ where $A,B,C$ are $\rats_p$-subspaces of $\octs_p$ of dimensions $i,j,k$ respectively, satisfying $N(A) = N(B) = N(C) = 0$ and $AB = BC = CA = 0$.  In other words the norm of any octonion in $A$ is $0$, and the product of any octonion in $A$ with any octonion in $B$ is $0$, etc...
\end{defn}
A number of remarks are in order, following work in \cite{Tit}.
\begin{itemize}
\item
The set of singular lines in $\octs_p$, e.g. ${\mathcal Fl}_{1,0,0}$, is a 6-dimensional quadric hypersurface in ${\mathbb P}^7$.
\item
The set of 3-dimensional hyperplanes in ${\mathcal Fl}_{1,0,0}$ (abr\'eg\'e in \cite{Tit}) come in two families, which may be identified with ${\mathcal Fl}_{0,1,0}$ and ${\mathcal Fl}_{0,0,1}$.  The symmetry between these two families and the points of the original quadric is known as triality.
\item
Incidence among points of ${\mathcal Fl}_{1,0,0}$ and points of ${\mathcal Fl}_{0,1,0}$ can be described as a point belonging to a 3-dimensional hyperplane, or as vanishing of octonionic multiplication.  This is described at the end of \cite{Tit}.
\item
Every two-dimensional subspace $A$ of $\octs_p$ satisfying $N(A)=0$ determines two other such subspaces, given by $B = \{ b \colon Ab = 0 \}$ and $C = \{ c \colon cA = 0 \}$.  Thus ${\mathcal Fl}_{2,0,0}$ is the same as ${\mathcal Fl}_{2,2,0}$ and ${\mathcal Fl}_{2,2,2}$.
\end{itemize}

The generalized flag varieties for $G$ can now be described via ${\mathcal Fl}_{i,j,k}$ for $i,j,k$ equal to $0,1,2$.  $G$ acts on ${\mathcal Fl}_{i,j,k}$ via the isotopy representation on $\octs_p^3$; we list the parabolic subgroups $P_{i,j,k}$ stabilizing a point on ${\mathcal Fl}_{i,j,k}$ for all $i,j,k$.  Note that the standard parabolic subgroups of $G$ are determined by subsets of the set $\{ \alpha_0, \ldots, \alpha_3 \}$ of simple roots.
\begin{center}
\begin{tabular}{|c|c|c|}
  \hline
  $i,j,k$ & Subset of $\{\alpha_0, \alpha_1, \alpha_2, \alpha_3 \}$ & Dimension of ${\mathcal Fl}_{i,j,k}$ \\
  \hline
  $0,0,0$ & $\{\alpha_0, \alpha_1, \alpha_2, \alpha_3 \}$ & $0$ \\
  $1,0,0$ & $\{\alpha_0, \alpha_2, \alpha_3 \}$ & $6$ \\
  $1,1,0$ & $\{\alpha_0, \alpha_3 \}$ & $9$ \\
  $1,1,1$ & $\{\alpha_0 \}$ & $11$ \\
  $2,1,1$ & $\{\alpha_1 \}$ & $11$ \\
  $2,2,1$ & $\{ \alpha_1, \alpha_2 \}$ & $10$ \\
  $2,2,2$ & $\{\alpha_1, \alpha_2, \alpha_3 \}$ & 9 \\
  \hline
 \end{tabular}
\end{center}

\subsection{The minimal representation of $\variety{E}_{7}(\rats_p)$}
We write $E' = \variety{E}_{7,3}(\rats_p)$, noting that $\variety{E}_{7,3}$ is split over $\rats_p$.  Let $Q_E = M_E N_E$ denote the (standard) maximal parabolic of $E'$ with abelian unipotent radical $N_E$.  $N_E$ can be identified with the exceptional Jordan algebra $J_3 = J_3(\octs_p)$.  Also, let $Q = MN$ denote the standard Borel subgroup of $SL_2^3$, consisting of triples of upper-triangular matrices;  $M$ is a split torus of rank $3$, and $N$ is unipotent abelian of dimension $3$.  There is a dual pair embedding $SL_2^3 \times_\nu G \hookrightarrow E'$ so that $Q_E \cap (SL_2^3 \times_\nu G) = Q \times_\nu G$.  The three-dimensional unipotent radical $N$ of $Q$ is identified with the diagonal elements of the Jordan algebra $J_3$.  This inclusion yields an orthogonal subspace $J_3^{\perp N}$ consisting of elements of the Jordan algebra $J_3$ with zeroes along the diagonal.  Thus $J_3^{\perp N}$ is identified with the set of triples of octonions $\octs_p^3$.

The Levi factor $M_E$ is isomorphic to $CE_6$, and acts via the 27-dimensional minuscule representation on $J_3$.  Let $\Omega'$ denote the orbit of a highest weight vector in $J_3$.  Define ${\Omega'}^{\perp N} = \Omega' \cap J_3^{\perp N}$.  $\Omega'$ is precisely the set of rank 1 elements of $J_3$, which implies:
\begin{lem}
The set ${\Omega'}^{\perp N}$ is the set of triples $(\alpha, \beta, \gamma)$ of octonions (not all of which are zero) such that $\alpha \beta = \beta \gamma = \gamma \alpha = 0$, and $N(\alpha) = N(\beta) = N(\gamma) = 0$.
\end{lem}

We can use this lemma to break ${\Omega'}^{\perp N}$ into locally closed pieces, based on the vanishing of $\alpha$, $\beta$, or $\gamma$.  Namely, for $i,j,k$ equal to $0$ or $1$, let $\mc{E}_{i,j,k}$ be the fibre bundle over $\mc{Fl}_{i,j,k}$ whose fibre over a point $(A,B,C)$ (a triple of subspaces of $\octs_p$ of dimensions $i,j,k$) is the set of triples $(\alpha, \beta, \gamma)$ spanning $(A,B,C)$.  Then we have the decomposition:
$${\Omega'}^{\perp N} = \bigsqcup_{i,j,k \in \{ 0,1 \} } \mc{E}_{i,j,k}.$$
Not all of the $i,j,k$ can equal zero, since the point $(0,0,0)$ is not in ${\Omega'}^{\perp N}$.

Let $\Pi'$ denote the minimal representation of $E'$.  By Theorem 1.1 of \cite{M-S}, the co-invariants of $\Pi'$ along the opposite unipotent radical $\bar N$ may be computed:
\begin{equation}
0 \rightarrow C_c^\infty({\Omega'}^{\perp N}) \rightarrow \Pi'_{\bar N} \rightarrow \Pi'_{\bar N_E} \rightarrow 0.
\label{SES7}
\end{equation}
Restricting $\Pi'$ to the dual pair $SL_2^3 \times_\nu G$, and taking co-invariants as above yields a representation of $M \times_\nu G = \Gp_m^3 \times_\nu G$.  The action of $\Gp_m^3 \times_\nu G$ on the terms in the above exact sequence is described by Theorem 1.1 of \cite{M-S}, and we recall this here.

Let $Isot$ denote the ``isotopy'' representation of $G$ on $\octs_p^3$, i.e., the direct sum of the three inequivalent 8-dimensional representations.  Let $\Gp_m^3$ also act on $\octs_p^3$ by scaling in the obvious way.  We write $Isot^1$ for the resulting representation of $\Gp_m^3 \times_\nu G$ on $\octs_p^3$.  The aforementioned results of Magaard and Savin \cite{M-S} imply:
\begin{prop}
The action of $\Gp_m^3 \times_\nu G$ on $C_c^\infty({\Omega'}^{\perp N})$ is given by:
$$[(t, g)f](\alpha, \beta, \gamma) = \vert t_1 t_2 t_3 \vert^{-4} f(Isot^1(t^{-1},g^{-1})(\alpha, \beta, \gamma)),$$
where $t = (t_1, t_2, t_3) \in \Gp_m^3$ and $g \in G$.
\end{prop}
The space $\Pi'_{\bar N_E}$ also admits a representation of $\Gp_m^3 \times_\nu G$, also described in Theorem 1.1 of \cite{M-S}.  The Levi component $M_E$ is isomorphic to $CE_6$, and we let $\Pi''$ denote the minimal representation of $E_6$, extended to $CE_6$ by having the center act trivially, and restricted to the pair $\Gp_m^3 \times_\nu G$ in $CE_6$.  Then we have:
\begin{prop}
The action of $\Gp_m^3 \times_\nu G$ on $\Pi'_{\bar N_E}$ is given by the representation:
$$\Pi'_{\bar N_E} \isom \left( \Pi'' \otimes \vert t_1 t_2 t_3 \vert^{-2} \right) \oplus \vert t_1 t_2 t_3 \vert^{-4}.$$
\end{prop}

\subsection{Tempered spherical representations}
Suppose that $\tau$ is a tempered spherical representation of $SL_2^3$, with regular parameter, so that:
$$\tau = \Ind_{\bar Q}^{SL_2^3} \chi_1(t_1) \chi_2(t_2) \chi_3(t_3) \vert t_1 t_2 t_3 \vert^{-1}.$$
The characters $\chi_i$ must be unitary and unramified.

We also consider irreducible smooth representations $\pi$ of $G$.  If $\pi$ is tempered spherical, there exist unramified characters of $\rats_p^\times$, $\psi_0, \psi_1, \psi_2, \psi_3$ so that $\pi$ occurs in $\Ind_{\bar B}^G \prod \psi_i \delta_{\bar B}^{1/2}$.
\begin{defn}
A tempered spherical representation $\pi$ is called {\it centrifugal} if it occurs in $\Ind_{\bar B}^G \prod \psi_i \delta^{1/2}$ with $\psi_0$ the trivial character.
\end{defn}
In other words, centrifugal representations are those whose Satake parameters are in the image of the inclusion of Langlands dual groups $PGL_2^3 \rightarrow {}^L G$.  For centrifugal representations, the induction from $T$ naturally factors through the three-step parabolic subgroup $P_{1,1,1}$; the Levi component $L_{1,1,1}$ of $P_{1,1,1}$ is the set of quadruples $(s, t_1, t_2, t_3)$ with $s \in GL_2$, $t_i \in \Gp_m$, and $det(s) \cdot t_1 t_2 t_3 = 1$.
\begin{prop}
Suppose that $\pi$ is centrifugal, with parameters $\psi_1, \psi_2, \psi_3$ ($\psi_0$ trivial by definition).  Then $\pi$ occurs in the induced representation:
$$Ind_{\bar P_{1,1,1}}^G \psi_1(t_1) \psi_2(t_2) \psi_3(t_3) \vert t_1 t_2 t_3 \vert^{-3}.$$
\end{prop}
\proof
This follows from inducing in stages, and an elementary computation of the modular character for $P_{1,1,1}$.
\qed

We refer to the triple $(\psi_1, \psi_2, \psi_3)$ as the parameter for a tempered spherical centrifugal representation.

Now, let $\pi$ be any irreducible smooth representation of $G$.  Applying Frobenius reciprocity, we have:
$$\Hom_{SL_2^3 \times_\nu G}(\Pi', \tau \boxtimes \pi) = \Hom_{\Gp_m^3 \times_\nu G} \left( \Pi'_{\bar N}, (\chi_1(t_1) \chi_2(t_2) \chi_3(t_3) \vert t_1 t_2 t_3 \vert^{-1}) \boxtimes \pi \right).$$
Since the characters $\chi_i$ are unitary, and $\vert t_1 t_2 t_3 \vert^{-1}$ is distinct from the two characters $\vert t_1 t_2 t_3 \vert^{-2}$ and $\vert t_1 t_2 t_3 \vert^{-4}$, we immediately get from the short exact sequence $\ref{SES7}$:
$$\Hom_{SL_2^3 \times_\nu G}(\Pi', \tau \boxtimes \pi) = \Hom_{\Gp_m^3 \times_\nu G} \left( C_c^\infty({\Omega'}^{\perp N}), (\chi_1(t_1) \chi_2(t_2) \chi_3(t_3) \vert t_1 t_2 t_3 \vert^{-1}) \boxtimes \pi \right).$$
Furthermore, we claim:
\begin{lem}
The only part of $C_c^\infty({\Omega'}^{\perp N})$ that contributes in the above equality is $C_c^\infty(\mc{E}_{1,1,1})$.  That is,
\begin{eqnarray*}
\Hom_{\Gp_m^3 \times_\nu G} \left( C_c^\infty({\Omega'}^{\perp N}), (\chi_1(t_1) \chi_2(t_2) \chi_3(t_3) \vert t_1 t_2 t_3 \vert^{-1}) \boxtimes \pi \right) \\
= \Hom_{\Gp_m^3 \times_\nu G} \left( C_c^\infty(\mc{E}_{111}), (\chi_1(t_1) \chi_2(t_2) \chi_3(t_3) \vert t_1 t_2 t_3 \vert^{-1}) \boxtimes \pi \right).
\end{eqnarray*}
\end{lem}
\proof
We look at the possible central characters of $C_c^\infty(\mc{E}_{i,j,k})$ for $i,j,k \in \{ 0, 1 \}$, viewed as representations of $\Gp_m^3 \times_\nu G$.  The possible characters are tabulated below (up to triality symmetry):
\begin{center}
\begin{tabular}{|c|c|}
  \hline
  $i,j,k$ & Central character \\
  \hline
  $1,0,0$ & $\rho_1(t_1) \vert t_1 t_2 t_3 \vert^{-4}$\\
  $1,1,0$ & $\rho_1(t_1) \rho_2(t_2) \vert t_1 t_2 t_3 \vert^{-4} $\\
  $1,1,1$ & $\rho_1(t_1) \rho_2(t_2) \rho_3(t_3) \vert t_1 t_2 t_3 \vert^{-4} $ \\
  \hline
 \end{tabular}
\end{center}
The characters $\rho_i$ in the above table may be arbitrary smooth characters of $\Gp_m$.  These sets are disjoint from the central character $(\chi_1(t_1) \chi_2(t_2) \chi_3(t_3) \vert t_1 t_2 t_3 \vert^{-1})$ when $\chi_i$ are unramified unitary, except when $i = j = k = 1$.  The lemma follows.
\qed

Using this lemma, it is not hard to prove a local theta correspondence:
\begin{thm}
Let $\tau$ be a tempered spherical representation of $SL_2^3$, with regular parameter.  Let $\pi$ be an irreducible spherical representation of $G$.  Then $\tau \boxtimes \pi$ occurs as a quotient of the restriction to $SL_2^3 \times_\nu G$ of $\Pi'$ if and only if $\pi$ is a tempered spherical centrifugal representation, whose parameters match those of $\tau$.
\end{thm}
\proof
From the lemma, we see that $\tau \boxtimes \pi$ occurs as a quotient of $\Pi'$ if and only if $((\chi_1(t_1) \chi_2(t_2) \chi_3(t_3) \vert t_1 t_2 t_3 \vert^{-1})) \boxtimes \pi$ occurs as a quotient of $C_c^\infty(\mc{E}_{1,1,1})$.  Now $\mc{E}_{1,1,1}$ is a fibre bundle over $\mc{Fl}_{1,1,1}$, with all fibres isomorphic to $\Gp_m^3$, and where $\Gp_m^3$ acts fibrewise in the obvious way, and the Levi component of $P_{1,1,1}$, $L_{1,1,1} \isom \Gp_m^3 \times SL_2$ acts fibrewise via the obvious representations of $\Gp_m^3 \subset L_{1,1,1}$.

Hence, we arrive at the following description of $C_c^\infty(\mc{E}_{1,1,1})$:
$$C_c^\infty(\mc{E}_{1,1,1}) \isom \vert t_1 t_2 t_3 \vert^{-4} \Ind_{\Gp_m^3 \times_\nu P_{1,1,1}}^{\Gp_m^3 \times_\nu G} C_c^\infty(\Gp_m) \boxtimes C_c^\infty(\Gp_m) \boxtimes C_c^\infty(\Gp_m).$$
The $C_c^\infty(\Gp_m)$ are essentially regular representations.  Thus $\tau \boxtimes \pi$ occurs as a quotient of the restriction of $\Pi'$ if and only if:
$$\Hom( \Ind_{P_{1,1,1}}^G \vert t_1 t_2 t_3 \vert^{-3} \chi_1(t_1) \chi_2(t_2) \chi_3(t_3), \pi) \neq 0.$$
Since we assume that $\pi$ is spherical, and the above induced representation has a unique spherical constituent, the theorem follows.
\qed
\begin{rem}
It is possible that the above theorem holds for general irreducible smooth representations $\pi$.  However, this seems to require an analysis of the reducibility of degenerate principal series induced from $P_{1,1,1}$.  As there are many (52, according to a MAPLE-assisted computation) $P_{1,1,1}$ double-cosets in $G$, this analysis seems cumbersome without any additional insight.
\end{rem}

 \subsection{The minimal representation of $\variety{E}_{8}(\rats_p)$}

The previous work on the minimal representation for $E' = \variety{E}_{7,3}(\rats_p)$ was a good warm-up to the more technical but similar work for $E = \variety{E}_{8,4}(\rats_p)$.  The techniques are essentially the same; however, instead of distinguishing representations by central character, we must instead use strategically chosen elements of the Bernstein center.

Noting that $\variety{E}_{8,4}$ is split over $\rats_p$, and let $P_E = L_E H_E$ denote the Heisenberg parabolic of $E$ as described for instance in \cite{G-W}.  The unipotent radical $H_E$ of $P_E$ has center $Z$, and $F = H_E / Z$ is a 56-dimensional $\rats_p$ vector space.  Let $P = LH$ denote the Heisenberg parabolic of $G$, so that $P = P_{2,2,2}$ in the previous.  $H$ is 9-dimensional, with one-dimensional center $Z$, and $H/Z$ is the vector space $C$ of 2 by 2 by 2 cubes over $\rats_p$.  There is a dual pair embedding $G \times_{\nu} G \hookrightarrow E$ so that $P_E \cap (G \times_{\nu} G) = P \times_{\nu} G$, and the centers $Z$ of $H_E$ and $H$ are identified.

The 56-dimensional space $F = H_E / Z$ can be viewed as the space of 2 by 2 matrices $\Matrix{x}{A_+}{A_-}{y}$ where $x,y \in \rats_p$, and $A_\pm$ are contained in the exceptional Jordan algebra $J_3(\octs_p)$.  The inclusion of the 8-dimensional space of cubes $C$ in $F$ yields an orthogonal subspace $F^{\perp C}$.  It is not hard to see that $F^{\perp C}$ consists of 2 by 2 matrices as above, where $x = y = 0$, and all diagonal entries of $A_\pm$ are $0$.  Hence elements of $F^{\perp C}$ are sextuples of octonions:
$$F^{\perp C} = \{ (\alpha_\pm, \beta_\pm, \gamma_\pm) \in \octs_p^6 \}.$$

The Levi factor $L_E$ of $P_E$ is isomorphic to $CE_7$, and acts via the 56-dimensional minuscule representation on $F$.  Let $\Omega$ denote the orbit of a highest weight vector in $F$.  Define $\Omega^{\perp C} = \Omega \cap F^{\perp C}$.  An analysis identical to that in Section 7 of \cite{M-S} yields:
\begin{lem}
The set $\Omega^{\perp C}$ is the set of sextuples $(\alpha_\pm, \beta_\pm, \gamma_\pm)$ of octonions such that if $A,B,C$ are the spans of $\alpha_\pm, \beta_\pm, \gamma_\pm$ respectively, then $N(A) = N(B) = N(C) = 0$ and $AB = BC = CA = 0$.  The sextuple $\alpha_\pm = \beta_\pm = \gamma_\pm = 0$ is excluded from $\Omega^{\perp C}$
\end{lem}

From this lemma, we may break up $\Omega^{\perp C}$ into a finite number of locally closed subsets.  For $i,j,k$ between 0 and 2, let $\mc{E}_{i,j,k}$ be the fibre bundle over $\mc{Fl}_{i,j,k}$ whose fibre over a point $(A,B,C)$ (a triple of subspaces of dimensions $i,j,k$ in $\octs_p$) is the set of sextuples $(\alpha_\pm, \beta_\pm, \gamma_\pm)$ spanning $(A,B,C)$ as in the lemma above.  We see immediately that:
$$\Omega^{\perp C} = \bigsqcup_{i,j,k \in \{ 0,1,2 \} } \mc{E}_{i,j,k}.$$
Again, not all of $i,j,k$ can equal $0$ in this decomposition.

Let $\Pi$ denote the minimal representation of $E$.  Recalling Theorem 6.1 of \cite{M-S}, the co-invariants of $\Pi$ along the unipotent subgroup $\bar Z$ opposite to $Z$ may be decomposed:
$$0 \rightarrow C_c^\infty(\Omega) \rightarrow \Pi_{\bar Z} \rightarrow \Pi_{\bar H_E} \rightarrow 0.$$
Furthermore, taking the co-invariants along all of $\bar H$, we get:
\begin{equation}
0 \rightarrow C_c^\infty(\Omega^{\perp C}) \rightarrow \Pi_{\bar H} \rightarrow \Pi_{\bar H_E} \rightarrow 0.
\label{SES}
\end{equation}

Restricting $\Pi$ to the dual pair $G \times_\nu G$, and taking co-invariants shows that the three terms in the above short exact sequence are representations of $L \times_\nu G$.  These representations of $L \times_\nu G$ are described by Theorem 6.1 of \cite{M-S}.  We begin with the action on $C_c^\infty(\Omega^{\perp C})$.

Recall $Isot$ is the ``isotopy'' representation of $G$ on the 24-dimensional space of triples of octonions $(\alpha, \beta, \gamma)$.  Let $St_{\octs_p}$ denote the ``standard'' action of $SL_2$ on pairs $(\kappa_+, \kappa_-)$ of octonions, given by:
$$\Matrix{a}{b}{c}{d} \left( {{\kappa_+} \atop {\kappa_-}} \right) = \left( {{a \kappa_+ + b \kappa_-} \atop {c \kappa_+ + d \kappa_-}} \right).$$
Let $Isot^2$ denote the resulting 48-dimensional representation of $SL_2^3 \times_\nu G$ on sextuples of octonions $(\alpha_\pm, \beta_\pm, \gamma_\pm)$.  In other words, $G$ acts on the triples $(\alpha_+, \beta_+, \gamma_+)$ and $(\alpha_-, \beta_-, \gamma_-)$ via the isotopy representation, and the three $SL_2$'s act on the three pairs $\alpha_\pm, \beta_\pm, \gamma_\pm$ via the standard representation.  Extend $Isot^2$ to $L \times_\nu G$ by letting the central $\Gp_m$ in $L$ act by uniformly scaling $\alpha_\pm, \beta_\pm, \gamma_\pm$.  Then Theorem 6.1 of \cite{M-S} implies:
\begin{prop}
The action of $L \times_\nu G$ on $C_c^\infty(\Omega^{\perp C})$ is given by:
$$[(l, g)f](\alpha_\pm, \beta_\pm, \gamma_\pm) = \vert \det \vert^{-5} f( Isot^2(l^{-1},g^{-1})(\alpha_\pm, \beta_\pm, \gamma_\pm)),$$
where $\det$ denotes the determinant character on $L$.
\label{Act}
\end{prop}

The space $\Pi_{\bar H_E}$ also admits a representation of $L \times_\nu G$, which is described in Theorem 6.1 of \cite{M-S}.  Let $\Pi'$ denote the minimal representation of $E_7$, extended to $CE_7$ by having the center act trivially, and restricted to the dual pair $L \times_\nu G$ in $CE_7$.  Again, $\det$ will denote the determinant character on $L$.  Then we have:
\begin{prop}
The action of $L \times_\nu G$ on $\Pi_{\bar H_E}$ is given by the representation:
$$\Pi_{\bar H_E} \isom \left( \Pi' \otimes \vert \det \vert^{-3} \right) \oplus \vert \det \vert^{-5}.$$
\label{Quo}
\end{prop}

\subsection{The Bernstein center for $L$}
We follow the methods of Section 4 of Magaard-Savin \cite{M-S} to pick out certain representations of $L$, using certain elements of the Bernstein center of $L$.  Recalling that $L$ is the group of triples $l = (l_1, l_2, l_3)$ of matrices in $GL_2$ such that $det(l_1) = det(l_2) = det(l_3)$, we begin by writing down a basis for the lattices $X_\bullet(T)$ and $X^\bullet(T)$.  As a basis for $X_\bullet(T)$, we choose:
\begin{eqnarray*}
\lambda_0(t) & = & \left( \Matrix{1}{0}{0}{t}, \Matrix{1}{0}{0}{t}, \Matrix{1}{0}{0}{t} \right), \\
\lambda_1(t) & = & \left( \Matrix{t}{0}{0}{t^{-1}}, \Matrix{1}{0}{0}{1}, \Matrix{1}{0}{0}{1} \right), \\
\lambda_2(t) & = & \left( \Matrix{1}{0}{0}{1}, \Matrix{t}{0}{0}{t^{-1}}, \Matrix{1}{0}{0}{1} \right), \\
\lambda_3(t) & = & \left( \Matrix{1}{0}{0}{1}, \Matrix{1}{0}{0}{1}, \Matrix{t}{0}{0}{t^{-1}} \right). \\
\end{eqnarray*}

We describe a basis for $X^\bullet(T)$ as follows:  if $l = (l_1, l_2, l_3) \in T$, then we let $\chi_0(l)$ denote the common determinant of $l_1, l_2, l_3$.  If $l_i = \Matrix{a_i}{0}{0}{b_i}$, then we define $\chi_i(l) = a_i$, for $i = 1,2,3$.  Then the canonical pairing $X_\bullet(T) \times X^\bullet(T) \rightarrow \ints$ satisfies:
\begin{eqnarray*}
\langle \lambda_i, \chi_i \rangle & = & 1, \mbox{ for } i = 0,1,2,3, \\
\langle \lambda_i, \chi_j \rangle & = & 0, \mbox{ for } i \neq j.
\end{eqnarray*}

The component of the Bernstein center acting non-trivially on representations generated by their Iwahori-fixed vectors is isomorphic to:
$$\complex[x_0, x_0^{-1}, x_1, x_1^{-1}, x_2, x_2^{-1}, x_3, x_3^{-1}]^W,$$
where $W$ is the abelian group of order 8 generated by $w_1, w_2, w_3$ acting on the $x_i$ by:
\begin{eqnarray*}
w_i x_i & = & x_i^{-1} x_0 \mbox{ for } i = 1,2,3, \\
w_i x_0 & = & x_0 \mbox{ for } i = 1,2,3, \\
w_i x_j & = & x_j \mbox{  for  } i \neq j, j \neq 0. \\
\end{eqnarray*}
Here the variables $x_i$ are identified with the cocharacters $\lambda_i$, but we use multiplicative notation for the variables $x_i$, rather than additive notation for $\lambda_i$.

If $E$ is a subquotient of an induced representation (from $T$ to $L$) with parameter $\chi = (\chi_0, \ldots, \chi_3)$, i.e., induced from the character $\chi$ extended to the Borel subgroup, then we have:
\begin{eqnarray*}
x_i + x_0 x_i^{-1} \vert_E & = & \chi_i(p) + \chi_0(p) \chi_i^{-1}(p), \\
x_0 \vert_E & = & \chi_0(p).
\end{eqnarray*}
If $\tau$ is a tempered spherical representation of $L$, then the parameter of $\tau$ has the form
$$\chi = (\chi_0 \vert \cdot \vert^{-3/2}, \chi_1 \vert \cdot \vert, \chi_2 \vert \cdot \vert, \chi_3 \vert \cdot \vert),$$
with all $\chi_i$ unitary characters.

Define the following elements of the Bernstein center of $L$:
\begin{eqnarray*}
T_i & = & x_i + x_0 x_i^{-1} \mbox{ for } i = 1,2,3, \\
T_0 & = & x_0.
\end{eqnarray*}
Then the $T_i$ act on the tempered $\tau$ above by:
\begin{eqnarray*}
T_i \vert_\tau & = & \chi_i(p) p^2 + \chi_0(p) \chi_i(p)^{-1} p^{-5/2}, \mbox{ for } i = 1,2,3, \\
T_0 \vert_\tau & = & \chi_0(p) p^{-3/2}.
\end{eqnarray*}

\subsection{Tempered spherical representations}

Suppose that $\pi$ is a irreducible tempered spherical representation of $G$ with regular parameter.  Then there exists a tempered spherical representation $\tau$ of $L$ so that:
$$\pi = \Ind_{\bar P}^G \left( \tau \otimes \vert \det \vert^{-5} \right) .$$

Let $\pi'$ be any other smooth representation of $G$.  Then by Frobenius reciprocity, we have:
$$\Hom(\Pi, \pi \boxtimes \pi') = \Hom(\Pi_{\bar H}, \left( \tau \otimes \vert \det \vert^{-5} \right) \boxtimes \pi').$$

Since $\tau$ is tempered, it has unitary central character.  Thus, there are no non-zero homomorphisms, or non-trivial extensions, from $\Pi' \otimes \vert \det \vert^{-3}$ to $\left( \tau \otimes \vert \det \vert^{-5} \right) \boxtimes \pi'$, since the central characters are disjoint.

We consider when the parameter for $\pi$ is regular, so that $\pi$ occurs as an irreducible induced representation.
\begin{prop}
Suppose that $\pi$ is a irreducible tempered spherical representation of $G$.  If we are given $\pi$ as an induced representation:
$$\pi = \Ind_{\bar P}^G \left( \tau \otimes \vert \det \vert^{-5} \right),$$
then we have:
$$\Hom(\Pi, \pi \boxtimes \pi') = \Hom(C_c^\infty(\Omega^{\perp C}), \left( \tau \otimes \vert \det \vert^{-5} \right) \boxtimes \pi').$$
\label{MDec}
\end{prop}
\proof
This is immediate from the short exact sequence \ref{SES} and Frobenius reciprocity.  Note that tempered representations are disjoint from the trivial representation, so that $\Hom(\pi, 1) = \Ext(\pi, 1) = 0$.
\qed

Finally, for $\pi$ induced from $\tau$ as before, $\tau$ must be tempered spherical, and the parameter of $\tau$ has the form
$$\chi = (\chi_0 \vert \cdot \vert^{-3/2}, \chi_1 \vert \cdot \vert, \chi_2 \vert \cdot \vert, \chi_3 \vert \cdot \vert),$$
with all $\chi_i$ unitary characters.

The action of the elements $p^5 T_i$ on $\tau \otimes \vert \det \vert^{-5}$ is given by:
\begin{eqnarray*}
p^5 T_i & = & \chi_i(p) p^2 + \chi_0(p) \chi_i(p)^{-1} p^{-5/2}, \\
p^5 T_0 & = & \chi_0(p) p^{-3/2}.
\end{eqnarray*}

\subsection{The p-adic correspondence}

The action of the elements $T_i$ on $C_c^\infty(\mc{E}_{i,j,k})$ can be explicitly computed, since $L$ acts fibrewise on the bundle $\mc{E}_{i,j,k}$ over $\mc{Fl}_{i,j,k}$.  If $\tau'$ is a spherical representation of $L$ occurring as a subquotient of $C_c^\infty(\mc{E}_{i,j,k})$, we denote its parameters by $(\rho_0, \ldots, \rho_3)$.  By knowing the possible parameters for $\tau'$, we tabulate the possible non-zero eigenvalues of $T_i$ on subquotients of $C_c^\infty(\mc{E}_{i,j,k})$ below, normalizing the $T_i$ by factors of $p^5$:

\begin{center}
\begin{tabular}{|c|c|c|}
  \hline
  $i,j,k$ & Eigenvalues of $p^5 T_0$ & Eigenvalues of $p^5 T_1, p^5 T_2, p^5 T_3$ \\
  \hline
  $1,0,0$ & $\rho_0(p)$ & $\rho_0(p) + 1,1 + \rho_0(p),1 + \rho_0(p)$ \\
  $1,1,0$ & $\rho_0(p)$ & $\rho_0(p) + 1$, $\rho_0(p) + 1$, $1 + \rho_0(p)$  \\
  $1,1,1$ & $\rho_0(p)$ & $\rho_0(p) + 1$, $\rho_0(p) + 1$, $\rho_0(p) + 1$ \\
  $2,1,1$ & $\rho_0(p)$ & $\rho_1(p) + {{\rho_0} \over {\rho_1}}(p)$, $\rho_0(p) + 1$, $\rho_0(p) + 1$  \\
  $2,2,1$ & $\rho_0(p)$ & $\rho_1(p) + {{\rho_0} \over {\rho_1}}(p)$, $\rho_2(p) + {{\rho_0} \over {\rho_2}}(p)$, $\rho_0(p) + 1$ \\
  $2,2,2$ & $\rho_0(p)$ & $\rho_1(p) + {{\rho_0} \over {\rho_1}}(p)$, $\rho_2(p) + {{\rho_0} \over {\rho_2}}(p)$, $\rho_3(p) + {{\rho_0} \over {\rho_3}}(p)$  \\
  \hline
 \end{tabular}
\end{center}

Looking at the above table, the part of every representation $C_c^\infty(\mc{E}_{i,j,k})$ generated by Iwahori-fixed vectors is an eigenspace for the operator $p^5(T_3 - T_0)$ (or a suitable variation under triality symmetry) of eigenvalue equal to $1$, except when $i = j = k = 2$.  On the other hand, if $\tau$ is an irreducible tempered spherical representation as before, then the eigenvalue of $p^5 (T_3 - T_0)$ equals:
$$v(\tau) = \chi_3(p) p^2 - \chi_0(p) p^{-3/2} + \chi_0(p) \chi_3(p)^{-1} p^{-5/2}.$$
Since the $\chi_i$ are unitary, an elementary estimate yields:
$$\vert v(\tau) \vert \geq p^2 - p^{-3/2} - p^{-5/2} \geq p^2 - 2 \geq 2.$$
Thus $v(\tau)$ cannot equal $1$.  Lemma 2.5 of \cite{M-S} now yields:
\begin{lem}
Let $\tau$ be an irreducible tempered spherical representation.  Then $\tau \otimes \vert \det \vert^{-5}$ occurs as a quotient of $C_c^\infty(\Omega^{\perp C})$ if and only if it occurs as a quotient of $C_c^\infty(\mc{E}_{2,2,2})$.
\end{lem}
It is now an easy step to get to the following:
\begin{thm}
Suppose that $\pi$ is an irreducible tempered spherical representation of $G$, whose parameter $(\chi_0, \ldots, \chi_3)$ is regular so that $\Ind_B^G \chi$ is irreducible.  If $\pi'$ is any irreducible smooth representation of $G$, then $\pi \boxtimes \pi'$ occurs as a quotient of the restriction to $G \times_\nu G$ of the minimal representation $\Pi$ if and only if $\pi \isom \pi'$.
\end{thm}
\proof
By Proposition \ref{MDec}, we have
$$\Hom(\Pi, \pi \boxtimes \pi') = \Hom(C_c^\infty(\Omega^{\perp C}), \left( \tau \otimes \vert \det \vert^{-5} \right) \boxtimes \pi').$$
From the last lemma, this yields:
$$\Hom(\Pi, \pi \boxtimes \pi') = \Hom(C_c^\infty(\mc{E}_{2,2,2}), \left( \tau \otimes \vert \det \vert^{-5} \right) \boxtimes \pi').$$
From Proposition \ref{Act}, we know that as a representation of $L \times_\nu G$,
$$C_c^\infty(\mc{E}_{2,2,2}) \isom \vert \det \vert^{-5} \otimes \Ind_{L \times_\nu P_{2,2,2}}^{L \times_\nu G} C_c^\infty(GL_2)^{\boxtimes 3} .$$
The representation $C_c^\infty(GL_2)^{\boxtimes 3}$ is essentially a regular representation of $L \times_\nu L$, and so every irreducible representation of $L \times_\nu L$ occurring as a quotient has the form $\tau \boxtimes \tau$ for some irreducible smooth $\tau$, and all such $\tau \boxtimes \tau$ occur.  Thus by Frobenius reciprocity again, we have:
$$\Hom(\Pi, \pi \boxtimes \pi') \neq 0 \mbox{ iff } \Hom(\tau \boxtimes \tau, \tau \boxtimes (\pi')_H) \neq 0.$$
By a result, attributed to Bernstein, proven by Bushnell in Theorem 3 of \cite{Bus}, we have:
$$\Hom(\tau, (\pi')_H) = \Hom(\Ind_{\bar P}^G \tau, \pi').$$
By the irreducibility of $\Ind_{\bar P}^G \tau = \pi$, we have
$$\Hom(\tau, (\pi')_H) = \Hom(\pi, \pi').$$
Finally, we get:
$$\Hom(\Pi, \pi \boxtimes \pi') \neq 0 \mbox{ iff } \Hom(\tau \boxtimes \pi, \tau \boxtimes \pi') = \Hom(\pi, \pi') \neq 0,$$
and the theorem follows.
\qed

It seems likely that some of the assumptions in this local theta correspondence could be removed with more technical work.  The assumption of regularity is necessary, since the $R$-group can be non-trivial for $G$ by a result of Keys \cite{Key}.  However, G. Savin has mentioned that working with the adjoint form would eliminate the $R$-group, and perhaps the need for the regularity assumption with it.  The spherical assumption could also likely be weakened, since in the regular representation of $L$ used above, various cuspidal representations occur paired with themselves.  We leave these details however, until a time when they might be necessary.

\section{Global theta correspondence}

In this section, we study global theta correspondences for the same groups we studied locally in the last section.  Namely, we hope to lift modular forms on $\variety{G}_c$ to modular forms on $\variety{G}_s$, and to holomorphic modular forms on $\variety{SL}_2^3$.

\subsection{The exceptional Jordan algebra}
Recall that $\Omega_c$ is the Coxeter's ring of integral octonions.  From $\Omega_c$, it is possible to define the exceptional Jordan algebra $J_3$ over $\ints$ (more precisely, the Jordan composition is defined over $\ints[1/2]$).  Let $\scheme{J}_3$ be the scheme over $\ints$ underlying the rank $27$ $\ints$-lattice of 3 by 3 Hermitian symmetric matrices over $\Omega_c$.  An element of $\scheme{J}_3(\ints)$ can be written in the form:
$$A = \left(
\begin{array}{ccc}
  a & \gamma & \bar \beta \\
  \bar \gamma & b & \alpha \\
  \beta & \bar \alpha & c \\
\end{array}%
\right)
,$$
where $a,b,c \in \ints$ and $\alpha, \beta, \gamma \in \Omega_c$.  $\scheme{J}_3$ naturally has additional structures.  Following \cite{GE1}, we define the cubic form $Det(A)$ by:
$$Det(A) = abc + Tr(\alpha \beta \gamma) - a \cdot N(\alpha) - b \cdot N(\beta) - c \cdot N(\gamma).$$
We define the adjoint matrix by:
$$A^\sharp = \left(
\begin{array}{ccc}
  bc - N(\alpha) & \bar \beta \bar \alpha - c \gamma & \gamma \alpha - b \bar \beta \\
  \alpha \beta - c \bar \gamma & ca - N(\beta) & \bar \gamma \bar \beta - a \alpha \\
  \bar \alpha \bar \gamma - b \beta & \beta \gamma - a \bar \alpha & ab - N(\gamma) \\
\end{array}%
\right)
.$$

If $R$ is a ring, then an element $A \in \scheme{J}_3(R)$ is said to have rank $1$ if $A \neq 0$, but its adjoint $A^\sharp = 0$.  There are no trace $0$ rank $1$ elements of $\variety{J}_3(\reals)$.

In \cite{Fre}, Freudenthal describes the $56$-dimensional representation of $\variety{E}_7(\reals)$ from $J_3$; following his construction, we define $\scheme{F}$ to be the group scheme over $\ints$ underlying the rank 56 $\ints$-lattice of $2$ by $2$ matrices of the form:
$$\phi = \left(
\begin{array}{cc}
  x & A_+ \\
  A_- & y \\
\end{array}
\right)
,$$
where $x,y \in \ints$ and $A_+, A_- \in \scheme{J}_3(\ints)$.

\subsection{Automorphic theta modules}
The automorphic realization of the global minimal representation $\Pi'$ of $\variety{E}_{7,3}(\adeles)$ follows from the work of Kim in \cite{Kim}.  Thus we have a map:
$$\Theta' \colon \bigotimes_v \Pi_v' \rightarrow L^2(\variety{E}_{7,3}(\rats) \backslash \variety{E}_{7,3}(\adeles)).$$
The representations $\Pi_v'$ are the minimal representations of $\variety{E}_{7,3}(\rats_v)$ for every place $v$ of $\rats$.

The work of Gan \cite{Ga2} gives an automorphic realization of the global minimal representation of $\variety{E}_{8,4}$.
$$\Theta \colon \bigotimes_v \Pi_v \rightarrow L^2(\variety{E}_{8,4}(\rats) \backslash \variety{E}_{8,4}(\adeles)).$$
The representations $\Pi_v$ are the minimal representations of $E_{8,4}(\rats_v)$ for every place $v$ of $\rats$.

The global minimal representations are spherical at every finite place (by \cite{Ga2} for $\variety{E}_{8,4}$), and the tensor product is taken with respect to suitably normalized spherical vectors.  The minimal $K$-type of the minimal representation of $\variety{E}_{7,3}(\reals)$ is one-dimensional, so in this case there is a natural way to choose a vector at the real place as well, up to scaling.  Both global minimal representations arise as quotients of globally parabolically induced representations from characters; thus there are natural {\it normalized} spherical vectors $t_p$, $t_p'$ for all finite primes $p$, for $\Pi_p, \Pi_p'$ respectively.  A vector $t \in \Pi$ (or $t' \in \Pi'$) is said to be standard if there is a decomposition $t = \bigotimes t_v$ (or $t' = \bigotimes t_v'$), where $t_v$ (or $t_v'$) is the normalized spherical vector for almost all $v$.

Fix $t$ a standard section of $\Pi$.  Let $t'$ denote the global normalized vector, spherical at all finite places, in $\Pi'$ considered by Kim \cite{Kim}.  Let $\theta = \Theta(t)$, and $\theta' = \Theta'(t')$.  Thus $\theta$ is an automorphic form on $\variety{E}_{8,4}$ and $\theta'$ is an automorphic form on $\variety{E}_{7,3}$.  We consider the Fourier expansion of $\theta$ and $\theta'$ along the following parabolic subgroups:  first, let $\variety{Q}_E = \variety{M}_E \variety{N}_E$ denote the (standard) maximal parabolic subgroup of $\variety{E}_{7,3}$ with abelian unipotent radical.  The unipotent radical $\variety{N}_E$ can naturally be identified with the exceptional Jordan algebra over $\rats$:  $\variety{N}_E \isom \variety{J}_3$.  Second, let $\variety{P}_E = \variety{L}_E \variety{H}_E$ denote the Heisenberg parabolic of $\variety{E}_{8,4}$.  The derived subgroup of $\variety{L}_E$ is the group $\variety{E}_{7,3}$.  The unipotent radical $\variety{H}_E$ is 57-dimensional, with 1 dimensional center $\variety{Z}$.  The 56-dimensional quotient can be identified with the Freudenthal space $\variety{F}$ described before; the action of the derived subgroup of $\variety{L}_E$ is the minuscule 56-dimensional representation discussed in \cite{Fre}.

Fix highest weight vectors $w'$ in $\variety{J}_3(\rats)$ and $w$ in $\variety{F}(\rats)$ for the action of the derived subgroups of $\variety{M}_E(\rats)$ and $\variety{L}_E(\rats)$ respectively.  Let $\Omega'$ and $\Omega$ denote the orbits of these highest weight vectors.  The orbit $\Omega'$ consists precisely of rank 1 elements of $\variety{J}_3(\rats)$.  The orbit $\Omega$ is more difficult to describe, but can be found in \cite{M-S},  Lemma 7.5.

We naturally identify $\variety{F}(\rats)$ with the set of characters of $\variety{H}_E(\adeles)$ which are trivial $\variety{H}_E(\rats)$.  Also, we identify $\variety{J}_3(\rats)$ with the set of characters of $\variety{J}_3(\adeles)$ which are trivial on $\variety{J}_3(\rats)$.

For a character $\phi \in \variety{J}_3(\rats)$, the Fourier coefficient of $\theta'$ is defined by:
$$\theta_\phi'(g) = \oint_{\variety{J}_3} \theta'(ng) \overline{\phi(n)} dn.$$

From the arguments in Section 5, Subsection 3, of \cite{G-S}, we have:
\begin{prop}
If $\phi$ is non-trivial, and $\theta'$ is non-zero, then $\theta_\phi'$ is non-zero only if $\phi$ lies in $\Omega'$ and equivalently has rank $1$.  If $\variety{M}_\phi$ denotes the stabilizer of $\phi$ in the Levi component $\variety{M}_E = \variety{E}_{6,2}$, then $\theta_\phi'(cg) = \theta_\phi'(g)$ for all $c \in \variety{M}_\phi(\adeles)$.
\end{prop}

When $t' = \bigotimes t_v'$ is the normalized spherical vector of $\Pi'$ at all finite places, and $t_\infty$ a well-chosen vector in the one-dimensional minimal $K$-type, Kim gives more precise information on the Fourier coefficients of $\theta'$ in \cite{Kim}:
\begin{prop}
The constant term of $\theta'$ is $1$, and the non-constant Fourier coefficients are given by non-zero constants (i.e., constant functions on $\variety{E}_{6,2}$) for all $\phi \in \scheme{J}_3(\ints)$:
$$a_\phi = 240 \sum_{d \vert c(\phi)} d^3,$$
where $c(\phi)$ is the largest integer such that $c(\phi)^{-1} \phi$ is in $\scheme{J}_3(\ints)$.
\end{prop}

Now we consider the Fourier expansion of $\theta$ along the Heisenberg parabolic $\variety{H}_E$, following Section 6 of \cite{Ga2}.  The $\variety{Z}$-constant term of $\theta$ is defined by:
$$\theta_{\variety{Z}}(g) = \oint_\variety{Z} \theta(zg) dz.$$
Suppose $\phi \in \variety{F}(\rats)$ is viewed as a character of $\variety{H}_E(\adeles)$ trivial on $\variety{H}_E(\rats)$.  Then the $\phi$-Fourier coefficient of $\theta$ is defined by
$$\theta_\phi(g) = \oint_\variety{F} \theta_{\variety{Z}}(ng) \overline{\phi(n)} dn.$$
From \cite{Ga2}, we know:
\begin{prop}
If $\phi$ is non-trivial, and $\theta$ is non-zero, then $\theta_\phi$ is non-zero only if $\phi$ is in the rational orbit $\Omega$ of a highest weight vector under the action of $\variety{L}_E(\rats)$.
\end{prop}
Furthermore, the constant term, when $\phi = 0$, is given by:
\begin{prop}
The constant term of $\theta$ along $\variety{H}_E$ is an automorphic form on the derived subgroup of $\variety{L}_E$ (which is isomorphic to $\variety{E}_{7,3}$);  $\theta_{\variety{H}_E}$ is given by $\theta_{\variety{H}_E} = c + \theta'$ where $c$ is a constant function and $\theta'$ is contained in the image of the automorphic realization of the minimal representation of $\variety{L}_E(\adeles)$.
\end{prop}

\subsection{Global theta lift to $\variety{SL}_2^3$}
We begin by considering the global theta lift from $\variety{G}_c$ to $\variety{SL}_2^3$.  Our methods are the same as those of \cite{Ga2}.  Fix a standard section $t'$ of $\Pi'$, and let $\theta' = \Theta'(t')$ as before.  Also fix an automorphic form $f_c$ on $\variety{G}_c$.  The theta lift of $f_c$ via $\theta'$ is defined by
$$\Phi(g) = \int_{\variety{G}_c} \theta'(g,g') f_c(g') dg',$$
where $g \in \variety{SL}_2^3(\adeles)$ and $(g, g')$ is considered as an element of $\variety{E}_{7,3}(\adeles)$ via the dual pair embedding:
$$\variety{SL}_2^3 \times_\nu \variety{G}_c \hookrightarrow \variety{E}_{7,3}.$$

The unipotent radical $\variety{N}$ of the Borel subgroup $\variety{Q}$ of $\variety{SL}_2^3$ is abelian, three-dimensional, and spanned by three subgroups $\variety{N}_{I}$, $\variety{N}_{II}$, $\variety{N}_{III}$.  For $\Phi$ to be cuspidal, we must examine the constant term of $\Phi$ along each of the three subgroups $\variety{N}_I, \variety{N}_{II}, \variety{N}_{III}$.  We have:
\begin{eqnarray*}
\Phi_{\variety{N}_I}(g) & = & \oint_{\variety{N}_I} \oint_{\variety{G}_c} \theta'(n g, g') f_c(g') dg' dn, \\
& = & \oint_{\variety{N}_I} \oint_{\variety{G}_c}  \sum_{\phi \in \variety{J}_3(\rats)} \theta_\phi'(n g , g') f_c(g') dg' dn,
\end{eqnarray*}
where the Fourier coefficients $\theta_\phi'$ satisfy $\theta_\phi'(n g, g') = \phi(n) \theta_\phi'(g, g')$ for all $n \in \variety{J}_3(\adeles)$.  Since $\variety{G}_c(\rats) \backslash \variety{G}_c(\adeles)$ compact, and the Fourier expansion of $\theta$ along $\variety{J}_3$ converges absolutely, the integration over $\variety{N}_I(\rats) \backslash \variety{N}_I(\adeles)$ may be brought inside.  This integration over $\variety{N}_I$ kills most of the Fourier coefficients of $\theta'$.  The remaining coefficients $\phi$ are those of the form:
$$\phi = \left(%
\begin{array}{ccc}
  0 & 0 & 0 \\
  0 & b & \alpha \\
  0 & \bar \alpha & c \\
 \end{array}%
\right),
$$
for $\alpha \in \octs_\rats$, $b,c \in \rats$, and $N(\alpha) = bc$.  We arrive at:
$$\Phi_{\variety{N}_I}(g) = \oint_{\variety{G}_c} f_c(g') \sum_{\phi \colon N(\alpha) = bc} \theta_\phi'(g, g') dg'.$$
If $\phi$ has the above form, with $N(\alpha) = bc = 0$, then the stabilizer in $\variety{G}_c$ of $\phi$ is all of $\variety{G}_c$.  Otherwise, the stabilizer in $\variety{G}_c$ of a $\phi$ of the form above is the $Spin_7$ subgroup $\variety{G}_I(\alpha)$.
The terms where $N(\alpha) = bc = 0$ may be expressed in terms of a period:
$$\sum_{\phi \colon \alpha = bc = 0} \oint_{\variety{G}_c} f_c(g') \theta_\phi'(g, g') dg'
= \sum_{\phi \colon \alpha = bc = 0} \theta_\phi'(g,1) \oint_{\variety{G}_c} f_c(g') dg'.$$
As long as $f_c$ is not $\variety{G}_c$-distinguished, i.e., $f_c$ is orthogonal to constant functions, the above integral vanishes.

For $N(\alpha) = bc \neq 0$, the stabilizer in $\variety{G}_c$ of $\phi$ is the group $\variety{G}_I(\alpha)$.  Moreover, $\variety{G}_c(\rats)$ acts transitively on the set of octonions $\alpha$ of a given non-zero norm $bc$.  Hence, assuming $f_c$ is not $\variety{G}_c$-distinguished, we can unfold the integral:
\begin{eqnarray*}
\Phi_{\variety{N}_I}(g) & = & \sum_{\phi \colon N(\alpha) = bc \neq 0} \oint_{\variety{G}_c} f_c(g') \theta_\phi'(g, g') dg' \\
& = & \sum_{bc \neq 0} \int_{\variety{G}_I(\alpha)(\adeles) \backslash \variety{G}_c(\adeles)} \theta_{\phi_{bc}}'(g, g') \oint_{\variety{G}_I(\alpha)} f_c(hg') dh dg'.
\end{eqnarray*}
Here $\phi_{bc}$ is a fixed $\phi$ satisfying $N(\alpha) = bc \neq 0$.  The inner integral is precisely the period of $f_c$ along the subgroup $\variety{G}_I(\alpha)$.  Hence we have shown:
\begin{prop}
If $f_c$ is not $\variety{G}_I$ (resp. $\variety{G}_{II}, \variety{G}_{III}$) distinguished, then the theta-lift $\Phi$ of $f_c$ is cuspidal along $\variety{N}_I$ (resp. $\variety{N}_{II}, \variety{N}_{III}$).
\end{prop}

Given the above sufficient condition for cuspidality of $\Phi$, we study its non-vanishing.  For generic $a,b,c \in \rats$, corresponding to a character $\psi = \psi_{a,b,c}$ of $\variety{N}(\adeles)$, the $\psi$-Fourier coefficient of $\Phi$ is given by:
$$\Phi_\psi(g) = \oint_{\variety{G}_c} f_c(g') \sum_{diag(\phi) = (a,b,c) } \theta_\phi'(g,g') dg',$$
where
$$\phi = \left(%
\begin{array}{ccc}
  a & \gamma & \bar \beta \\
  \bar \gamma & b & \alpha \\
  \beta & \bar \alpha & c \\
\end{array}%
\right).$$
If the diagonal entries of $\phi$ are $(a,b,c)$ as above, so $\phi$ restricts to $\psi$ on $\variety{N}$, we write $Res(\phi) = \psi$.  For $\phi$ to have rank $1$ (these are the only characters for which $\theta'$ has non-vanishing coefficients), we must have $N(\alpha) = bc, N(\beta) = ca, N(\gamma) = ab$ and $\alpha \beta = c \bar \gamma$ as well.  Thus the stabilizer of such $\phi$ in $\variety{G}_c$ is the $G_2$ subgroup $\variety{G}_2(\alpha, \beta, \gamma)$ of $\variety{G}_c$.

\begin{lem}
The group $\variety{G}_c(\rats)$ acts transitively on the set of $\phi$ satisfying $Res(\phi) = \psi$ as above.
\label{TrL}
\end{lem}
\proof
Our proof directly follows suggestions of the referee:  the $\variety{G}_c(\rats)$-orbits on these $\phi$ correspond to elements of the kernel of the canonical map:
$$H^1(Gal(\bar \rats / \rats), \variety{G}_2(\alpha, \beta, \gamma)(\bar \rats)) \rightarrow H^1( Gal(\bar \rats / \rats), \variety{G}_c(\bar \rats)).$$
Since both of the groups $\variety{G}_2(\alpha, \beta, \gamma)$ and $\variety{G}_c$ are simply-connected, the first Galois cohomology vanishes over every finite place.  Applying the Hasse principle, we see that the orbits are classified by the kernel of the map:
$$H^1(Gal(\complex / \reals), \variety{G}_2(\alpha, \beta, \gamma)(\complex)) \rightarrow H^1(Gal(\complex / \reals), \variety{G}_c(\complex)).$$
For anisotropic groups $G$ over $\reals$, the first Galois cohomology yields the set of conjugacy classes of elements of $G(\reals)$ of order 1 or 2.  The preimage of the trivial conjugacy class in $\variety{G}_c(\reals)$ is again the trivial conjugacy class, i.e., the distinguished element, of $\variety{G}_2(\alpha, \beta, \gamma)(\reals)$.  Hence the kernel of the above map is trivial.
\qed

Now unfolding the expression for $\Phi_\psi$ yields:
$$
\Phi_\psi(g) = \int_{\variety{G}_2(\alpha, \beta, \gamma)(\adeles) \backslash \variety{G}_c(\adeles)} \theta_{\phi_0}'(g, g') \oint_{\variety{G}_2(\alpha,\beta,\gamma)} f_c(h g') dh dg'.
$$
Here, $\phi_0$ is a fixed $\phi$ restricting to $\psi$.  The inner integral is again a period, from which we derive:
\begin{prop}
If $f_c$ is not $\variety{G}_2$-distinguished then $\Phi$ vanishes.  Inversely, if $f_c$ is $\variety{G}_2$-distinguished, then $\Phi$ does not vanish.
\label{G2P}
\end{prop}

The vanishing of $\variety{G}_2$-periods clearly implies the vanishing of $\Phi$ by the equation above.  We are left to check the inverse statement; suppose that $f_c$ is $\variety{G}_2(\alpha, \beta, \gamma)$-distinguished.  Note that:
$$\Phi_\psi(1) = \int_{\variety{G}_2(\adeles) \backslash \variety{G}_c(\adeles)} \theta_\phi'(g') \Period_{f_c}^{\variety{G}_2}(g') dg'.$$
The same analysis as in Section 5, Proposition 4.5 of \cite{G-S} can now be used to show that the above quantity does not vanish for a suitable choice of $\theta'$ in the image of $\Theta'$.

\subsection{Global theta lift to $\variety{G}_s$}

We can now study the theta lift from $\variety{G}_c$ to $\variety{G}_s$.  Fix $\theta$ in the image of $\Theta$, and an automorphic form $f_c \in \mc{A}_c = \mc{A}(\variety{G}_c)$.  The theta lift of $f_c$ via $\theta$ is the automorphic form on $\variety{G}_s$ defined by:
$$f_s(g) = \oint_{\variety{G}_c} f_c(g') \theta(g, g') dg',$$
where $g \in \variety{G}_s(\adeles)$, and where we view $(g, g')$ as an element of $\variety{E}_{8,4}(\adeles)$ via the dual pair embedding:
$$\variety{G}_s \times_\nu \variety{G}_c \hookrightarrow \variety{E}_{8,4}.$$

$f_s$ is a cusp form on $\variety{G}_s$ if the constant terms of $f_s$ along the unipotent radicals of all (standard) maximal parabolic subgroups of $\variety{G}_s$ vanish.  The maximal parabolic subgroups of $\variety{G}_s$ are the Heisenberg parabolic $\variety{P}$, and three parabolic subgroups $\variety{Q}_i = \variety{M}_i \variety{N}_i$ with $\variety{N}_i$ abelian of dimension $6$.

The description of the constant term $\theta_{\variety{H}_E} = c + \theta'$ for $\theta'$ in the image of $\Theta'$ shows that if $f_c$ is orthogonal to the constant functions, then:
\begin{eqnarray*}
(f_s)_\variety{H}(g) & = & \oint_{\variety{H}} \oint_{\variety{G}_c} f_c(g') \theta(h g,  g') dg' dh, \\
& = & \oint_{\variety{G}_c} f_c(g') \theta_{\variety{H}_E}(g, g') dg', \\
& = & \oint_{\variety{G}_c} f_c(g') \theta'(g, g') dg'.
\end{eqnarray*}
To deduce the second line from the first above, note that integration over $\variety{H}(\rats) \backslash \variety{H}(\adeles)$ kills all non-constant Fourier coefficients of $\theta$ except those $\phi \in \variety{F}(\rats)$ of the form $\psi = \Matrix{0}{A_+}{A_-}{0}$, with $A_\pm$ having zeroes along their diagonal.  But by Lemma 2.7 of \cite{Ga2}, and the fact that there are no trace zero, rank one elements of $\variety{J}_3(\rats)$, such $\psi$ do not exist; the only term to survive is the constant term along all of $\variety{H}_E$.

Coupled with the results in Proposition \ref{G2P}, we have now shown:
\begin{prop}
$(f_s)_\variety{H}$, as an automorphic form on $\variety{SL}_2^3$, is the theta lift of $f_c$ via $\theta'$. The constant term of the theta lift $f_s$ along the Heisenberg parabolic $\variety{H}$ vanishes if $f_c$ is not $\variety{G}_2$-distinguished.
\end{prop}

To check whether the lift $f_s$ is cuspidal, we must determine when $(f_s)_{\variety{N}_i}$ vanishes for the three abelian unipotent radicals $\variety{N}_i$.  We may write $\variety{N}_i = (\variety{N}_i \cap \variety{H}_E ) \oplus \variety{N}_i'$, for a suitable subgroup $\variety{N}_i'$.  Applying this decomposition, we compute these constant terms now:
\begin{eqnarray*}
(f_s)_{\variety{N}_i}(g)  & = & \oint_{\variety{N}_i} \oint_{\variety{G}_c} f_c(g') \theta(n g, g') dg' dn, \\
& = & \oint_{\variety{N}_i'} \oint_{\variety{N}_i \cap \variety{H}_E} \oint_{\variety{G}_c} f_c(g') \theta(n' n g,  g') dg' dn' dn, \\
& = & \oint_{\variety{N}_i'} \oint_{\variety{G}_c} f_c(g') \sum_{\psi \in \Omega_i} \theta_\psi(n' g, g') dg' dn'.
\end{eqnarray*}
In the above, $\Omega_i$ denotes the subset of the orbit $\Omega$ orthogonal to $\variety{N}_i \cap \variety{H}_E$.  Using Gan's classification of elements of $\Omega$ in Lemma 2.7 of \cite{Ga2}, an element of $\Omega_i$, for $i = 1$ looks like:  $\psi = \Matrix{0}{A_+}{A_-}{d}$ with $d \in \rats$ and:
$$A_+ = \left(%
\begin{array}{ccc}
  a & 0 & 0 \\
  0 & 0 & 0 \\
  0 & 0 & 0 \\
 \end{array}%
\right),
A_- = \left(%
\begin{array}{ccc}
  0 & 0 & 0 \\
  0 & b & \alpha \\
  0 & \bar \alpha & c \\
 \end{array}%
\right),$$
where $N(\alpha) = bc$.  The stabilizer of such an element in $\Omega_i$ is the $Spin_7$ subgroup $\variety{G}_I(\alpha)$.  From this we see:
\begin{prop}
If $f_c$ is not $\variety{G}_I$ (resp. $\variety{G}_{II}, \variety{G}_{III}$) distinguished, the constant terms $(f_s)_{N_1}$ (resp. $(f_s)_{N_2}$, $(f_s)_{N_3}$) vanish.
\end{prop}

If $f_c$ is not $\variety{G}_2$-distinguished, then it cannot be $\variety{G}_I$, $\variety{G}_{II}$, or $\variety{G}_{III}$ distinguished, since any of the latter subgroups contains a $\variety{G}_2$ subgroup.  Hence we see:
\begin{thm}
If $f_c$ is not $\variety{G}_2$-distinguished, then $f_s$ is cuspidal.
\end{thm}

Finally, we derive a condition for the theta lift to be non-vanishing.  For a generic cube $c_{ijk}$, corresponding to a character $\psi$ of $\variety{H}(\adeles)$, we compute the Fourier coefficient of $f_s$ at $\psi$.  Such a coefficient is given by:
$$(f_s)_\psi(g) = \oint_{\variety{G}_c} f_c(g') \sum_{Res(\phi) = \psi} \theta_\phi(g, g') dg'.$$
Those $\phi \in \Omega$ that restrict to a generic $\psi$ have the form $\phi = x \cdot \Matrix{1}{ A}{ A^\sharp}{ Det(A)}$, for $x \in \rats^\times$, $A \in \variety{J}_3(\rats)$, using the description of $\Omega$ in Lemma 2.7 of \cite{Ga2} again.  If $A \in \variety{J}_3(\rats)$ has the form:
$$A = \left( \begin{array}{ccc}
  a & \gamma & \bar \beta \\
  \bar \gamma & b & \alpha \\
  \beta & \bar \alpha & c \\
 \end{array}%
\right),$$
then the subgroup of $\variety{G}_c$ stabilizing $A$ is $\variety{SU}_3(\alpha, \beta, \gamma)$ (generically).  An easy computation shows that $\variety{SU}_3(\alpha, \beta, \gamma)$ stabilizes $A^\sharp$ as well.  Thus the stabilizer of $\phi$ is precisely $\variety{SU}_3(\alpha, \beta, \gamma)$.

As in Lemma \ref{TrL}, we can see that $\variety{G}_c(\rats)$ acts transitively on the set of $\phi$ restricting to $\psi$ as above.  For the orbits are given by the kernel in Galois cohomology of the map:
$$H^1(Gal(\bar \rats / \rats), \variety{SU}_3(\alpha, \beta, \gamma)(\bar \rats)) \rightarrow H^1(Gal(\bar \rats / \rats), \variety{G}_c(\bar \rats)).$$
Once again, both groups are simply connected, and anisotropic over $\reals$, so the kernel is trivial and there is one $\variety{G}_c(\rats)$ orbit.

We continue the computation of the Fourier coefficient by unfolding the integral:
$$(f_s)_\psi(g) = \int_{\variety{SU}_3(\alpha, \beta, \gamma)(\adeles) \backslash \variety{G}_c(\adeles)} \theta_{\phi_0}(g, g') \oint_{\variety{SU}_3(\alpha,\beta,\gamma)} f_c(h g') dh dg'.$$
Here $\phi_0$ is a fixed character restricting to $\psi$, and again the inner integral is a period:
$$(f_s)_\psi(g) = \int_{\variety{SU}_3(\alpha, \beta, \gamma)(\adeles) \backslash \variety{G}_c(\adeles)} \theta_{\phi_0}(g, g') \Period^{\variety{SU}_3(\alpha,\beta,\gamma)} f_c(g') dg'.$$

The same argument applies to degenerate Fourier coefficients, replacing the group $\variety{SU}_3(\alpha, \beta, \gamma)$, by a $\variety{G}_2$, $\variety{G}_I$, $\variety{G}_{II}$, or $\variety{G}_{III}$ subgroup.  In particular, the vanishing of $\variety{SU}_3(\alpha, \beta, \gamma)$ periods implies the vanishing of periods for such larger subgroups.  

Hence we have:
\begin{prop}
If $f_c$ is not $\variety{SU}_3$-distinguished, then all Fourier coefficients of the theta-lift $f_s$ of $f_c$, with respect to $\variety{H} / \variety{Z}$ vanish.
\end{prop}

The following lemma is adapted from a lemma in Section 8 of \cite{GGS}:

\begin{lem}
An automorphic form $f_s$ on $\variety{G}_s$ vanishes if and only if the constant term $(f_s)_{\variety{Z}}$ vanishes.
\end{lem}
\proof
The argument from Lemma 9.1 in \cite{Ga2} applies.  Namely, let $\variety{Z}_{1,1,1}$ denote the center of the unipotent radical of the three-step parabolic subgroup $\variety{P}_{1,1,1} \subset \variety{G}_s$.  Then $\variety{Z}_{1,1,1}$ is two-dimensional, and contains $\variety{Z}$.  Thus if $(f_s)_{\variety{Z}} = 0$, then for any character $\phi$ on $\variety{Z}_{1,1,1}$ trivial on $\variety{Z}$, we have $(f_s)_\phi = 0$.  But the Levi component $\variety{L}_{1,1,1}$ of $\variety{P}_{1,1,1}$ can be used to take any character $\phi$ of $\variety{Z}_{1,1,1}$ to a character trivial on $\variety{Z}$.  Thus every Fourier coefficient of $f_s$ along $\variety{Z}_{1,1,1}$ vanishes, and so $f_s$ vanishes.
\qed

We now have:  
\begin{thm}
The theta lift $f_s$ vanishes if and only if $f_c$ is not $\variety{SU}_3$-distinguished.
\end{thm}
\proof
If $f_c$ is not $\variety{SU}_3$-distinguished, then all Fourier coefficients of $f_s$ with respect to $\variety{H} / \variety{Z}$ vanish.  Hence $(f_s)_\variety{Z}$ vanishes, and by the previous lemma, $f_s$ vanishes.  The other direction follows from the same analysis as in \cite{G-S}.
\qed

\section{Examples of modular forms}

In this section, we discuss the construction of some modular forms on $\variety{G}_c$, and the resulting modular forms on $\variety{SL}_2^3$ and $\variety{G}_s$.  In particular, a non-zero constant function on $\variety{G}_c$ yields a modular form, whose theta-lift to $\variety{SL}_2^3$ and $\variety{G}_s$ should be identified with Eisenstein series.  We consider the significance of the Fourier coefficients of the lift to $\variety{G}_s$, achieving the arithmetic part of a Siegel-Weil formula.

\subsection{Spherical harmonics and invariants}

A useful way to construct modular forms on $\variety{G}_c$ is through the theory of spherical harmonics, and invariant polynomials.  Identify $\octs_c$ with the representation previously denoted $V_{0, \omega}$, $\omega = (1,0,0)$ of $\variety{G}_c(\reals)$.  Let $P(n, \octs_c)$ denote the space of homogeneous polynomials of degree $n$ on $\octs_c \simeq \reals^8$ with real coefficients.  We see that $P(n, \octs_c)$ is a representation of $Spin_8 = \variety{G}_c(\reals)$, identified with $Sym^n(V_{0, (1,0,0)})$.  Let $r^2$ denote the homogeneous polynomial of degree $2$ on $\octs_c$, given by the quadratic norm form.  Let $\triangle$ denote the Laplacian associated to the norm form, normalized so that $\triangle r^2 = r^2$.  

Define $H(n, \octs_c)$ to be the subspace of $P(n, \octs_c)$ consisting of homogeneous polynomials $p$ which are harmonic, i.e., satisfy $\triangle p = 0$.  Then the decomposition of $P(n, \octs_c)$ as a representation of $Spin_8$ is well-known:
$$P(n, \octs_c) = \bigoplus_{0 \leq m \leq \lfloor n/2 \rfloor} r^{2m} H(n-2m, \octs_c).$$
Each space of harmonic polynomials $H(n, \octs_c)$ is isomorphic as a $Spin_8$ representation to $V_{0, \omega}$ with $\omega = (n,0,0)$.

From Proposition \ref{MF1} it follows that:
\begin{prop}
The space of modular forms on $\variety{G}_c$ of level $1$, and of weight $(0, (n,0,0))$ can be identified with $H(n, \octs_c)^{\Gamma_c}$, i.e., with harmonic polynomials of degree $n$, invariant under $\scheme{G}_c(\ints)$.
\end{prop}

In order to study such invariant polynomials for the group $\Gamma_c$, we describe a close relationship between $\Gamma_c$ and the Weyl group $W_E$ of the $E_8$ root system.  First, we note that the finite group $\Gamma_c$ is a central extension of $\scheme{G}_c(\FF_2)$ by the abelian group of order $4$, $\nu(\ints)$:
$$1 \rightarrow \nu(\ints) \rightarrow \Gamma_c \rightarrow \scheme{G}_c(\FF_2) \rightarrow 1.$$

By our construction of $\scheme{G}_c$, we can consider the three images, $\Gamma_c^I, \Gamma_c^{II}, \Gamma_c^{III}$ of $\Gamma_c$ in $\scheme{SO}(\Omega_c, N)$.  Each one of these groups is still a central extension of $\scheme{G}_c(\FF_2)$, this time by a group of order $2$.  In particular, the action of $\Gamma_c$ on $\octs_c$ factors through the quotient $\Gamma_c^I$ and:
$$H(n, \octs_c)^{\Gamma_c} = H(n, \octs_c)^{\Gamma_c^I}.$$

By identifying the $E_8$ root lattice with Coxeter's octonions (scaling if necessary), the Weyl group $W_E$ acts on $\Omega_c$ by the reflection representation $ref$.  Let $W_E^+$ denote the kernel of $det \circ ref$ in $W_E$, the subgroup of index $2$ acting by proper isometries on $\Omega_c$.  It is known that $W_E^+$ is also a central extension:
$$1 \rightarrow \{ \pm 1 \} \rightarrow W_E^+ \rightarrow \scheme{G}_c(\FF_2) \rightarrow 1.$$
In fact, we have:
\begin{prop}
There is an isomorphism between $W_E^+$ and $\Gamma_c^I$ which intertwines the reflection representation of $W_E^+$ and the representation of $\Gamma_c^I$ on $\Omega_c$.  
\end{prop}

From this, we immediately get:
\begin{cor}
The space of modular forms on $\variety{G}_c$ of level $1$, and of weight $(0, (n,0,0))$ can be identified with $H(n, \octs_c)^{W_E^+}$.
\end{cor}

The ring of invariants for the full $E_8$ Weyl group $W_E$ is generated by invariants of degrees $2,8,12,14,18,20,24,30$.  The invariants for the index $2$ subgroup $W_E^+$ consist precisely of the invariants for $W_E$ and the ``skew invariants'' as described in \cite{Kan}.  The skew invariants form a free cyclic module over the ring of invariants, generated by a single skew invariant $Sk_{240}$ of degree $240$ (the number of reflections in $W_E$).  

One may choose canonical fundamental invariants $I_2 = r^2, I_8, \ldots, I_{30}$, up to scalar multiple, as described in \cite{Iwa}.  These canonical invariants will be harmonic, except for $I_2$, which has Laplace eigenvalue $1$. The skew invariant $Sk_{240}$ may be chosen to be harmonic as well, as discussed in \cite{Kan}.  Therefore, we see:
\begin{thm}
There are canonical, up to scalar multiple, modular forms $F_d$ on $\variety{G}_c$ of level $1$ and of weights $(0, (d,0,0))$ for $d \in \{8, 12, 14, 18, 20, 24, 30, 240 \}$, associated to the canonical harmonic invariant polynomials, $I_d$ and the skew invariant polynomial $Sk_{240}$.  The invariant $I_2$ of degree $2$ corresponds to the trivial representation of $Spin_8$, and yields a constant modular form.
\end{thm}

We may apply the results of Corollary \ref{PeC} to see that all of the modular forms $F_d$ are $\variety{G}_2$-distinguished.  Therefore, the theta-lifts of $F_d$ to $\variety{SL}_2^3$ and to $\variety{G}_s$ do not vanish.  In order to understand these $F_d$, and their theta-lifts, one must first be able to explicitly write down the canonical polynomial invariants $I_d$, preferably in a way that exploits the octonionic structure of the $E_8$ root lattice.  We leave this study to a future paper.
  
\subsection{Lifting the trivial modular form to $SL_2^3$}

Though the constant modular form on $\variety{G}_c$ is uninteresting by itself, its theta-lifts are worthy of study, especially considering the exceptional Siegel-Weil formula of Gan \cite{Ga1}.  We describe the theta-lifts to $\variety{SL}_2^3$ and $\variety{G}_s$ here, focusing on connections to Eisenstein series, and a description of Fourier coefficients.

The lifting of the trivial modular form on $\variety{G}_c$ to $\variety{SL}_2$ is particularly simple, given Kim's thorough description of a theta function on $\variety{E}_{7,3}$, as well as some more general work of Gross-Elkies in \cite{GE2}.  In particular, we can work classically throughout, beginning with the theta function on the exceptional tube domain:
$$\theta'(Z) = 1 + 240 \sum_{A \geq 0, rk(A) = 1} \left( \sum_{d \vert c(A)} d^3 \right) e^{2 \pi i \langle A, Z \rangle}.$$
Some explanation is necessary for the above formula.  We view $\theta'$ as a holomorphic function on the exceptional tube domain:
$${\mathcal D} = \{Z = X + iY \colon X \in \variety{J}_3(\reals), Y \in \variety{J}_3(\reals)_+ \}.$$
The summation is over elements $A \in \scheme{J}_3(\ints)$, of rank $1$, which are positive semi-definite.  If $A \in \scheme{J}_3(\ints)$, the integer $c(A)$ refers to the largest positive integer dividing $A$.  
The embedding of $\variety{G}_c \times_\nu \variety{SL}_2^3$ in $\variety{E}_{7,3}$, and the action of $\variety{E}_{7,3}(\reals)$ on the exceptional tube domain, allows us to define the theta-lift of the constant modular form on $\variety{G}_c$ as:
$$\Phi(z_1, z_2, z_3) = \int_{\variety{G}_c(\reals)} \theta' \left( g 
\left(%
\begin{array}{ccc}
  z_1 & 0 & 0 \\
  0 & z_2 & 0 \\
  0 & 0 & z_3 \\
\end{array}%
\right)
\right) 
dg.$$
For $Z$ a diagonal element of the exceptional tube domain, and $g \in \variety{G}_c(\reals)$ as in the above formula, $gZ = Z$.  Therefore, we see that:
$$e^{2 \pi i \langle A, g Z \rangle} = e^{2 \pi i \langle A, Z \rangle}.$$
Hence Kim's exceptional form is invariant under translation by $\variety{G}_c(\reals)$.  It follows immediately that the modular form $\Phi$ on $\variety{SL}_2^3$ is given by:
$$\Phi(z_1, z_2, z_3) = \sum_{a_1, a_2, a_3 \in \nats} \rho(a) e^{2 \pi i (a_1 z_1 + a_2 z_2 + a_3 z_3)},$$
where the coefficients $\rho$ are given by:
$$\rho(a) = 240 \sum_{rk(A) = 1, diag(A) = a} \left( \sum_{d \vert c(A)} d^3 \right) .$$
From Proposition 4.1 of \cite{GE2}, or noticing the connection between these coefficients and the theta function for the $E_8$ root lattice, we get:
\begin{thm}
The theta-lift of the constant function on $\variety{G}_c$ via Kim's exceptional form on $\variety{E}_{7,3}$ to a modular form on $\variety{SL}_2^3$ is the product of three Eisenstein series of weight $4$ for $\variety{SL}_2$:
$$\Phi(z_1, z_2, z_3) = E_4(z_1) E_4(z_2) E_4(z_3).$$
\end{thm}

This may be seen as an analogue of the Siegel-Weil formula for the dual pair $\variety{G}_c \times_\nu \variety{SL}_2^3$ in $\variety{E}_{7,3}$, though it does not contain significant new arithmetic information.

\subsection{Lifting the trivial modular form to $G_s$}
Far more interesting to us is the lifting of the trivial modular form to a modular form on $\variety{G}_s$.  For this, we fix $t_f = \bigotimes t_p$, the product of the normalized spherical vectors of the minimal representation of $\variety{E}_{8,4}(\rats_p)$ for all (finite) primes $p$.  By the archimedean theta correspondence of Loke, which we discussed in Theorem \ref{LoT}, we have a (unique up to scaling) embedding $\iota$ of the quaternionic discrete series representation $\pi_{10}$ of $\variety{G}_s(\reals)$, paired with the trivial representation of $\variety{G}_c(\reals)$ into the local minimal representation $\Pi_\infty$ of $\variety{E}_{8,4}(\reals)$.  

For any vector $v$ of the quaternionic discrete series $\pi_{10}$, embedded as a vector $\iota(v)$ in $\Pi_\infty$, we have an automorphic form $\Theta(t_f \otimes \iota(v))$ on $\variety{E}_{8,4}(\adeles)$.  Restricting from $\variety{E}_{8,4}(\adeles)$ to $\variety{G}_s(\adeles)$ yields an automorphic form $\Psi_v$ on $\variety{G}_s(\adeles)$, and the map $v \mapsto \Psi_v$ is a modular form of weight $10$ and level $1$, in the sense of Definition \ref{MoD}.  

Hereafter, we denote by $\Psi$ the modular form obtained by lifting the trivial modular form on $\variety{G}_c$ as above.  Since $\Psi$ has level one, there are well-defined Fourier coefficients $a_c$ for $\Psi$, indexed by 2 by 2 by 2 cubes $c$.  For simplicity, we consider only those coefficients corresponding to projective, non-degenerate cubes; in this case, the $SL_2(\ints)^3$-invariance of the Fourier coefficients allows us to consider only cubes in normal form as well.  

Following the work in Section 11 of \cite{Ga1}, we have, after suitably normalizing all coefficients:
\begin{prop}
If $c$ is a projective non-degenerate cube, then the Fourier coefficient $a_c$ is equal to the number of elements $\omega \in \Omega \cap \scheme{F}(\ints)$ which restrict to the cube $c$.
\end{prop}

The counting problem in this proposition is tractable, when $c$ is in normal form.  We let $c$ be the cube:
$$
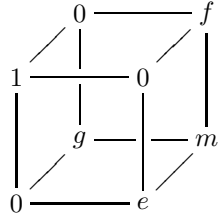
\xymatrix@!0{
  & 0 \ar@{ }[rr] \ar@{ }'[d][dd]
      &  & f \ar@{ }[dd]        \\
  1 \ar@{ }[ur]\ar@{ }[rr]\ar@{ }[dd]  
      &  & 0 \ar@{ }[ur]\ar@{ }[dd] \\
  & g \ar@{ }'[r][rr]
      &  & m                \\
  0 \ar@{ }[rr]\ar@{ }[ur]
      &  & e \ar@{ }[ur]        }
$$

 An element $\omega \in \Omega \cap \scheme{F}(\ints)$ restricting to $c$ is then a matrix $\Matrix{1}{A}{A^\sharp}{Det(A)}$, where the diagonal entries of $A^\sharp$ are $e,f,g$, $Det(A) = m$, and $A$ is given by:
$$A = \left(%
\begin{array}{ccc}
  0 & \gamma & \bar \beta \\
  \bar \gamma & 0 & \alpha \\
  \beta & \bar \alpha & 0 \\
\end{array}%
\right)
.$$
Hence we see that:
\begin{prop}
The number of elements $\omega \in \Omega \cap \scheme{F}(\ints)$ restricting to $c$ (and hence the Fourier coefficient $a_c$) is equal to the number of triples $(\alpha, \beta, \gamma) \in \Omega_c^3$ such that $N(\alpha) = -e$, $N(\beta) = -f$, $N(\gamma) = -g$, and $Tr(\alpha \beta \gamma) = m$.
\label{Em1}
\end{prop}

\subsection{An embedding problem}
The previous proposition gives some interpretation of the Fourier coefficients of the theta lift of the constant function on $\variety{G}_c$.  However, in analogy to classical Siegel-Weil formulas, and the exceptional Siegel-Weil formula of Gan \cite{Ga1}, we look for a more arithmetically interesting interpretation.  For this, we introduce the following algebraic object:
\begin{defn}
A $QT$-structure (over $\ints$) of rank $n$ consists of a free $\ints$-module $\Lambda$ of rank $n$, three integer-valued quadratic forms $Q_1, Q_2, Q_3$ on $\Lambda$, and a trilinear form $T \colon \Lambda \otimes \Lambda \otimes \Lambda \rightarrow \ints$.  
\end{defn}

We have already seen one $QT$-structure, namely Coxeter's integral octonions $\Omega_c$, endowed with the trilinear form $Tr(\alpha \beta \gamma)$, letting all three quadratic forms $Q_i$ be the norm form.  We can deduce from results of Bhargava in \cite{Bha}, that every 2 by 2 by 2 cube also yields a $QT$-structure as well.  We describe this construction here:

Suppose that the cube $c$ is in normal form as before, with discriminant $D \neq 0$, and let $R(D)$ denote the quadratic ring of discriminant $D$.  Associated to the cube $c$, we get three invertible oriented ideal classes $I_1, I_2, I_3$.  Let $\Lambda$ denote the free $\ints$-module of rank $2$, with basis $\lambda, \mu$.  By results in the Appendix of \cite{Bha}, there exist $\ints$-module isomorphisms from $\Lambda$ to $I_1, I_2, I_3$ (choosing appropriate representatives for these ideal classes), such that such that the quadratic norm form on $R(D) \otimes \rats$, applied to the images of $x \lambda + y \mu$ in $I_1, I_2, I_3$ yields the following binary quadratic forms:
\begin{eqnarray*}
Q_1 & = & -e x^2 + m xy + fg y^2, \\
Q_2 & = & -f x^2 + m xy + eg y^2, \\
Q_3 & = & -g x^2 + m xy + ef y^2.
\end{eqnarray*}
In this way, we get three quadratic forms $Q_1, Q_2, Q_3$ on the lattice $\Lambda$, by looking at the norm form on $I_1, I_2, I_3$.  In order to get a trilinear form, we look at the trilinear map which multiplies elements of $I_1, I_2, I_3$ (viewing them as elements of $R(D) \otimes \rats$), and applies the trace map from $R(D)$ to $\ints$.  With the same basis $x,y$ for $\ints^2$ as above, we may compute the values of the trilinear form; for example $T(\lambda,\lambda,\lambda) = m$ and $T(\mu,\mu,\mu) = m^2 + 2 efg$.  In any case, we see that a 2 by 2 by 2 cube, in normal form, yields a $QT$-structure $(\Lambda, Q, T)$ of rank $2$.  Moreover, with respect to a well-chosen basis $\lambda,\mu$ for this structure, we have $Q_1(\lambda) = -e$, $Q_2(\lambda) = -f$, $Q_3(\lambda) = -g$, and $T(\lambda,\lambda,\lambda) = m$.  With a few tedious arithmetic computations, and Proposition \ref{Em1}, we get:
\begin{thm}
The Fourier coefficients $a_c$ associated to any projective non-degenerate cube $c$ in normal form count the number of embeddings of the $QT$-structure associated to $c$ into the $QT$-structure coming from Coxeter's integral octonions.
\end{thm}
This can be seen as the arithmetic part of a Siegel-Weil formula for the dual pair $(\variety{G}_c, \variety{G}_s)$ in $\variety{E}_{8,4}$.  It would be interesting if one could identify the modular form $\Psi$ with an Eisenstein series on $\variety{G}_s$, and use this to obtain other formulae for the $a_c$, especially relating to values of L-functions.
\bibliography{D4Modular}
\end{document}